\documentclass[twocolumn]{autart}    % Enable this line and disable the 
\usepackage{graphicx}          % Include this line if your 
\newcommand{\col}[1]{\operatorname{col}(#1)}

\renewcommand{\d}{\mathrm{d}}
\renewcommand{\t}{^{\top}}
\renewcommand{\e}{\mathrm{e}}

\newcommand{\tint}{\textstyle\int}

\usepackage{pgfplots}
\pgfplotsset{compat=newest} % Kern
\newlength\figureheight % Kern
\newlength\figurewidth % Kern

\usepackage{environ}

\usepackage{psfrag}
\usepackage{amssymb}
\usepackage{graphicx}
\usepackage{amsmath,amsfonts,amsxtra}
\usepackage{mathrsfs}
\usepackage{mathtools}

\usepackage{tikz}
\usetikzlibrary{arrows,automata}
\allowdisplaybreaks[1]

\begin{document}

\begin{frontmatter}
%\runtitle{Insert a suggested running title}  % Running title for regular 
                                              % papers but only if the title  
                                              % is over 5 words. Running title 
                                              % is not shown in output.

\title{Cooperative output regulation for a network of \\ parabolic systems with varying  parameters} % Title, preferably not more 
                                                % than 10 words.

%\thanks[footnoteinfo]{This paper was not presented at any IFAC 
%meeting. Corresponding author M.~T.~Cicero. Tel. +XXXIX-VI-mmmxxi. 
%Fax +XXXIX-VI-mmmxxv.}

\author{Joachim Deutscher}\ead{joachim.deutscher@uni-ulm.de}% ~and~   % Add the Joachim Deutscher}\ead{joachim.deutscher@rt.eei.uni-erlangen.de
%\author{Jakob Gabriel}\ead{jakob.gabriel@fau.de}               % e-mail address 
%\author[Baiae]{Publius Maro Vergilius}\ead{vergilius@culture.ir}  % (ead) as shown

\address{Institut f\"ur Mess-, Regel- und Mikrotechnik, Universit\"at Ulm, 
	Albert-Einstein-Allee 41, D--89081 Ulm, Germany}  % Please supply                                              
%\address[Rome]{Senate House, Rome}             % full addresses
%\address[Baiae]{The White House, Baiae}        % here.
          
                             % keyword list or with the 
                                          % help of the Automatica 
                                          % keyword wizard

\begin{abstract}
This paper is concerned with the cooperative output regulation problem for a network of agents with different dynamics described by parabolic PDEs subject to spatially- and time-varying parameters. Firstly, a networked controller is designed achieving output synchronization for identical finite-dimensional reference models, which deliver the state of the global reference model required for the synchronization to the parabolic agents. The latter can be subject to local disturbances acting in-domain, on all boundaries and on the anti-collocated output to be controlled. The cooperative output regulation problem is solved by designing local output feedback regulators for the parabolic agents. This requires the solution of time-varying regulator equations and the design of disturbance observers for parabolic systems with spatially- and time-varying coefficients. For this, a systematic backstepping approach is provided and it is shown that cooperative output regulation with exponential convergence is ensured for the resulting multi-agent system. The results of the paper are applied to the cooperative output regulation of a heterogeneous network of four parabolic agents in the presence of local disturbances.
\end{abstract}

\end{frontmatter}

%###########################################################
\section{Introduction}
In a lot of applications spatially separated systems have to cooperatively solve a common control task by transmitting information through a communication network. In order to avoid an excessively large information exchange, the communication is limited to neighbouring systems leading to a restricted communication topology. This setup results in the \emph{networked control of multi-agent systems (MAS)}, which is still a very active research topic (see the recent monographs \cite{Bu19,Lu19} for MAS with finite-dimensional agents).\\
Since many applications require to take both the temporal and spatial system dynamics into account, it is also of interest to design networked controllers for MAS with  distributed-parameter agents. Applications include industrial furnaces consisting of a network of heaters (see  \cite{Ca14}), networks of HVAC systems in building climate control (see, e.\:g., \cite{Bo09} for distributed-parameter models and \cite{Sa17} for the corresponding MAS approach), networks of Lithium-Ion cells in battery management (see, e.\:g., \cite{Qu16} and \cite{Ta17} for distributed-parameter models)  or consensus control in environmental applications (see \cite{Tr12}). Another interesting topic  for applying distributed-parameter MAS is the formation control for networks of tricopters carrying a payload using heavy ropes (see, e.\:g., \cite{Ir19}). These are still emerging technologies, which require a strong methodical background for the networked control. Different from the finite-dimensional case, however, less results can be found in the literature for the control of distributed-parameter MAS. The recent contributions \cite{Dem13,Pi16} consider parabolic agents, while parabolic PDEs with a diffusive coupling are tackled in \cite{Wu12,Wu16}. For wave equations results can be found in \cite{Ag20,Che20}. General classes of distributed-parameter agents were dealt with in \cite{Dem18} by making use of an abstract setting. Common to these works is the fact that only homogeneous MAS are considered, which limits their applicability. In particular, no results can be found in the literature to synchronize heterogeneous MAS with distributed-parameter agents. This is of practical importance, because the agent dynamics may differ due to the production process or different environment. In addition, also heterogeneous agents may result from the modelling of unequal technological processes. It should be noted that recently PDE theory has also been applied to the networked control of finite-dimensional agents by making use of distributed-parameter continuum models (see, e.\:g., \cite{Fri11,Qi15,Fre20}). More specifically, the backstepping approach is utilized in these contributions for the deployment of the agents (see, e.\:g., \cite{Kr08} for the backstepping approach to PDEs). \textcolor{blue}{Different from these works not a continuum of agents is modelled in this paper by a single PDE, but the dynamics of each individual agent are described by a PDE.} %Therein, in particular, the work \cite{Fre20} considers continuum models with time and spatially varying parameters.\\
In order to provide a systematic approach for the networked control of heterogeneous MAS with distributed-parameter agents, this paper extends the \emph{cooperative output regulation problem} to infinite dimensions (see, e.\:g., \cite{Li19} for the finite-dimensional cooperative output regulation problem).  A general framework is obtained by taking parabolic agents into account, which differ both in their spatial and temporal parameters. By extending the two-step approach of \cite{Wi11,Su12a,Su12} to the distributed-parameter MAS in question, a cooperative output feedback regulator is systematically determined. In a first step, the design of a homogeneous network of local reference models, the so-called \emph{cooperative reference observer}, is presented. This result is obtained by adapting the corresponding methods in \cite{Su12a,Su12}. The reference input for all agents is specified by the solution of a global reference model. By making use of the communication network the local reference models synchronize with the global reference model so that they can deliver the state of the global reference model to the agents.  
%This directly generalizes the corresponding approach in \cite{Deu15a}, where only one ODE-leader (reference model) and one parabolic agent with spatially-varying coefficients (plant) are taken into account. 
%An advantage of considering finite-dimensional reference agents lies in the fact that only their lumped-parameter outputs have to be communicated through the network. This significantly reduces the communication load for the synchronization of the distributed-parameter agents. 
With this, it is possible to design local regulators for anti-collocated outputs in a second step so that each agent can follow the output of the global reference model. In addition, local disturbances are taken into account, that may act in-domain, on both boundaries and on the anti-collocated output to be controlled. These exogenous signals are described by local finite-dimensional disturbance models. This leads to a challenging output regulation problem for parabolic PDEs with spatially- and time-varying coefficients, which was not considered in the literature so far (see, e.\:g., \cite{Aul16} for the output regulation of distributed-parameter systems). Therefore, this result is also of interest on its own for parabolic systems with varying parameters.  %Related results for ODEs in \cite{Ich06,Zh06} consider both general and periodically time-varying systems. 
Related results for distributed-parameter systems are limited to  time periodic systems (see \cite{Pa17}). %Besides the networked control of MAS this output regulation problem is also of particular interest for applications, where the plant is described by a parabolic PDE with varying parameters. Such systems occur, e.\:g., due to the linearisation of nonlinear PDEs along desired trajectories (see \cite{Me07}). Moreover, PDEs with time-varying spatial domains may be transformed into PDEs with constant domain but space and time dependent coefficients (see \cite{Iz15}). 
The solution of the output regulation problem in question leads to regulator equations in form of a Cauchy problem for parabolic PDEs. This is different from the time-invariant case, where only boundary value problems for ODEs have to be solved (see \cite{Deu15a}). A systematic solution procedure for the Cauchy problem is provided by making use of the results in \cite{Ru08}. This determines the feedforward regulator ensuring output regulation. The corresponding feedback part for stabilizing the closed-loop system follows from applying the backstepping approach developed in \cite{Sm05,Vaz08,Me09,Me13}. With this, the local output feedback regulator is obtained by designing a disturbance observer for each agent and its disturbance model. Consequently, a backstepping observer has to be determined for a coupled parabolic ODE-PDE system, in which the PDE is subject to spatially- and time-varying coefficients. This new challenging problem is solved by extending the results in \cite{Deu15a} to this more general setup with the aid of the backstepping methods in \cite{Sm05,Vaz08,Me09,Me13}. In order to facilitate the design, the mapping into a stable PDE-ODE cascade is achieved by decoupling the PDE subsystem in backstepping coordinates. This yields a simple initial boundary value problem with constants coefficients for the second transformation step, while the first step only needs the usual kernel equations. Output synchronization with the global reference model is verified for the resulting networked controlled MAS with exponential convergence. This provides a general and scalable framework to solve cooperative output regulation problems for both homogeneous and heterogeneous networks of parabolic agents.\\
The next section introduces the considered cooperative output regulation problem. Then, the design of the cooperative reference observer is considered in Section \ref{sec:refnw}. The two subsequent sections present the local regulator design and the resulting cooperative output feedback regulator is investigated in Section \ref{sec:coufeedreg}. The results of the paper are illustrated for a heterogeneous network of four parabolic agents in the presence of local disturbances.

\emph{Elements from graph theory and definitions.} The communication topology between the agents is described by a \emph{digraph} $\mathcal{G} = \{\mathcal{V},\mathcal{E},A_{\mathcal{G}}\}$, in which $\mathcal{V}$ is a set of $N$ \emph{nodes} $\mathcal{V} = \{v_1,\ldots,v_N\}$, one for each agent and $\mathcal{E} \subset \mathcal{V} \times \mathcal{V}$ is a set of \emph{edges} that models the information flow from the node $v_j$ to $v_i$ by $(v_i,v_j) \in \mathcal{E}$. This flow is weighted by $a_{ij} \geq 0$, which are the element of the \emph{adjacency matrix} $A_{\mathcal{G}} = [a_{ij}] \in \mathbb{R}^{N \times N}$ with $a_{ii} = 0$, $i = 1,\ldots,N$. From this, the \emph{Laplacian matrix} $L_{\mathcal{G}}  \in \mathbb{R}^{N \times N}$ of the graph $\mathcal{G}$ can be derived by $L_{\mathcal{G}} = D_{\mathcal{G}} - A_{\mathcal{G}}$, where $D_{\mathcal{G}} = \operatorname{diag}(d_1,\ldots,d_N)$ with $d_i = \sum_{k=1}^Na_{ik}$, $i = 1,\ldots,N$. %It is assumed that there are no \emph{self-loops}, i.\:e., $a_{ii} = 0$, $i = 1,\ldots,N$, holds for $A_{\mathcal{G}}$. 
%Hence, the elements $l_{ij}$ of $L_{\mathcal{G}}$ are $l_{ij} = \sum_{k=1}^{N}a_{ik},$ $i = j$ and $l_{ij} =-a_{ij}$, $i \neq j$. 
A \emph{path} from the node $v_j$ to the node $v_i$ is a sequence of $r \geq 2$ distinct nodes $\{v_{l_1},\ldots,v_{l_r}\}$ with $v_{l_1} = v_j$ and $v_{l_r} = v_i$ such that $(v_k,v_{k+1}) \in \mathcal{E}$. A graph $\mathcal{G}$ is said to be \emph{connected} if there is a node $v$, called the \emph{root}, such that, for any node $v_i \in \mathcal{V} \setminus \{v\}$, there is a path from $v$ to $v_i$. For further details on graph theory see, e.\:g., \cite[Ch. 2]{Mes10}. Define the sets $\mathbb{R}^+_{t_0} = \{t \in \mathbb{R} \;|\; t > t_0\}$, $t_0 \geq 0$, and $\overline{\mathbb{R}}^{+}_{t_0} = \{t \in \mathbb{R} \;|\; t \geq t_0\}$. The function $f: \mathbb{R}_{t_0}^+  \to \mathbb{R}$ belongs to the \emph{Gevrey class $G_{\alpha}(\mathbb{R}^+_{t_0})$ of order $\alpha$} if $f \in C^{\infty}(\mathbb{R}^+_{t_0})$ and there exist positive constants $M$ and $R$ such that $
\sup_{t \in K}|f^{(i)}(t)| \leq \frac{M}{R^i}(i!)^{\alpha}, \quad \forall i \in \mathbb{N}_0$, 
for any compact subset $K \subset \mathbb{R}_{t_0}^+$ (see, e.\:g., \cite{Rod93}).

%###########################################################
\section{Problem formulation}\label{sec:probform}
Consider a \emph{multi-agent system (MAS)} consisting of the $N > 1$ \emph{heterogeneous  agents} described by the parabolic systems
\begin{subequations}\label{drs}
    \begin{align}
	 \!\!\dot{x}_i(z,t) &= \textcolor{blue}{\lambda_i(z)x''_i(z,t) + \phi_i(z)x'_i(z,t)} +  \! a_i(z,t)x_i(z,t)\nonumber\\
	                & \quad   + g_{1,i}\t(z,t)d_i(t), \quad (z,t) \in \textcolor{blue}{\overset{\circ}{\Omega}_i} \times \mathbb{R}^+_{t_0}\label{pdes}\\
          x'_i(0,t) \! &= \!q_{i}(t)x_i(0,t)  \! + \!  g_{2,i}\t(t)d_i(t), &&\hspace{-1.85cm} t \!\in \!\mathbb{R}^+_{t_0}\label{bc1}\\
	    \!\!\!\!\!\!x'_i(\ell_i,t) \! &= \textcolor{blue}{\!q^{\ell}_i(t)x_i(\ell_i,t)} \! +  \! \textcolor{blue}{b_i}u_i(t)  \! +  \! g_{3,i}\t(t)d_i(t),&&\hspace{-1.85cm}  t \!\in\! \mathbb{R}^+_{t_0} \label{bc2}
\end{align}
with the spatial domain \textcolor{blue}{$\Omega_i = [0,\ell_i]$}, $b_i \in \mathbb{R}  \setminus \{0\}$ and the \emph{output to be controlled} $y_i(t) \in \mathbb{R}$, which has not to be available for the controller, and the \emph{collocated measurement} $\eta_i(t) \in \mathbb{R}$ given by
\begin{align}
	           y_i(t) &= \textcolor{blue}{c_i}x_i(0,t) + g_{4,i}\t(t)d_i(t),&&  t \in \overline{\mathbb{R}}^+_{t_0}\label{outcontr}\\
     	    \eta_i(t) &= \textcolor{blue}{c_{m,i}}x_i(\ell_i,t),&&  t \in \overline{\mathbb{R}}^+_{t_0}\label{outmeas}	    
	\end{align}
\end{subequations}	
for $i = 1,\ldots,N$ and $c_i, c_{m_i} \in \mathbb{R} \setminus \{0\}$. The state is $x_i(z,t) \in \mathbb{R}$, \textcolor{blue}{$\lambda_i \in C^2(\Omega_i)$ satisfies $0 < \underline{\lambda}_i \leq \lambda_i(z) \leq \overline{\lambda}_i < \infty$ on $z \in \Omega_i$,  $\phi_i \in C^1(\Omega_i)$, $a_i \in C(\Omega_i) \times G_{\alpha}(\mathbb{R}^+_{t_0})$}, $1 \leq \alpha < 2$, are assumed. The known input locations of the  disturbance $d_i(t) \in \mathbb{R}^{m_i}$ are characterized by $g_{1,i} \in (C_p(\textcolor{blue}{\Omega_i}) \times G_{\alpha}(\mathbb{R}^+_{t_0}))^{m_i}$ with $C_p(\textcolor{blue}{\Omega_i})$ denoting the space of piecewise continuous functions on $\Omega_i$ and $g_{k,i} \in (G_{\alpha}(\mathbb{R}^+_{t_0}))^{m_i}$, $k = 2,3,4$. In \eqref{bc1} and \eqref{bc2} the coefficients $q_{i}, q^{\ell}_{i} \in G_{\alpha}(\mathbb{R}^+_{t_0})$ specify Robin or Neumann BCs, the input is $u_i(t) \in \mathbb{R}$ and the initial condition (IC) of the system reads $x_i(z,t_0) = x_{0,i}(z) \in \mathbb{R}$. \textcolor{blue}{In order to simplify the presentation of the results, the spatial and Hopf-Cole-type transformations given in \cite[Ch. 4.8]{Kr08} are used to normalize \eqref{drs} with $\lambda_i(z) = \lambda_i = const.$, $\phi_i = 0$ and $\ell_i = 1$, $i = 1,\ldots,N$. It can be verified that the networked control of the agents in the original coordinates also solves the posed cooperative output regulation problem.}\\
%\begin{rem}
%Note that by using the spatial and time transformations given in \cite[Ch. 7.1.1]{Me13} the more general linear parabolic PDE $\dot{x}_i(z,t) = (\delta_i(z)x'_i(z,t))' + \beta_i(z)x'_i(z,t) +  \gamma_i(z,t)x_i(z,t)+ \bar{g}_{1,i}(z,t)d_i(t)$ on the domain $(z,t) \in (0,L) \times \mathbb{R}^+$ with $\delta_i \in C^2[0,L]$, $\beta_i \in C^1[0,L]$, $\gamma_i \in C[0,L] \times G_{\alpha}(\mathbb{R}^+_{t_0})$ and $\bar{g}_{1,i} \in (C_p[0,L] \times G_{\alpha}(\mathbb{R}^+_{t_0}))^{m_i}$, $1 \leq \alpha < 2$, can be mapped to \eqref{pdes}. This yields the same BCs and outputs as in \eqref{bc1}--\eqref{outmeas}. Hence, the following results are also valid for this general system class, that describes many heat transfer problems or chemical processes (see \cite{Me13}).
%\end{rem}
For all agents a common \emph{reference input} $r(t) \in \mathbb{R}$ is specified by the solution of the known \emph{global reference model}
\begin{subequations}\label{smodelg}
\begin{align}
 \dot{w}(t) &= Sw(t), && t \in \mathbb{R}^+_{t_0}, \quad w(t_0) = w_0 \in \mathbb{R}^{n_w} \label{sigmodssg}\\
       r(t) &= p\t w(t), &&  t \in \overline{\mathbb{R}}^+_{t_0}\label{sigmoy_routg}
\end{align}
\end{subequations}
with $p \in \mathbb{R}^{n_w}$ and the pair $(p\t,S)$ observable.  It is assumed that the spectrum $\sigma(S)$ of $S \in \mathbb{R}^{n_w \times n_w}$ has only eigenvalues on the imaginary axis, i.\:e., $\sigma(S) \subset \text{j}\mathbb{R}$. Hence, \eqref{smodelg} describes a wide class of reference inputs including polynomial as well as trigonometric functions of time. The \emph{disturbance} $d_i$, $i = 1,\ldots,N$, acting on the individual agents is described by the known \emph{local disturbance model}
\begin{subequations}\label{smodel}
	\begin{align}
	\dot{v}_i(t) &= S_iv_i(t), &&\hspace{-0.2cm} t \in \mathbb{R}^+_{t_0}, \quad v_i(t_0) = v_{0,i} \in \mathbb{R}^{n_{v_i}} \label{sigmodss}\\
	d_i(t) &= P_iv_i(t), &&\hspace{-0.2cm}  t \in \overline{\mathbb{R}}^+_{t_0}\label{sigmoy_rout}
	\end{align}
\end{subequations}
with $P_i \in \mathbb{R}^{m_i \times n_{v_i}}$, the pair $(P_i,S_i)$ observable and the diagonalizable matrix $S_i \in \mathbb{R}^{n_{v_i} \times n_{v_i}}$ satisfying $\sigma(S_i) \subset \text{j}\mathbb{R}$.\\
The agents consist of two groups. The first group is composed of the agents $i$, $i = 1,\ldots,n$, $n \geq 1$, which have access to the reference input $r$ and are therefore called the \emph{informed agents}. In contrast, the information about the reference input can only be broadcast to the remaining agents $i$, $i = n+1,\ldots,N$, through a communication network with the digraph $\mathcal{G}$. %(see, e.\:g., \cite{Su12} for a quick reference on graph theory and \cite{Bu19,Lu19} for further results). 
In particular, these agents have only access to the reference information of their neighbours due to the communication constraints. Hence, they are the so-called \emph{uninformed agents}. As a consequence, a stabilizing cooperative regulator is required, in order to ensure output regulation, i.\:e.,
%output regulation.
%In this paper, the \emph{cooperative output regulation problem} is solved by utilizing  observer-based feedforward control. This amounts to designing a stabilizing cooperative output feedback regulator such that
\begin{equation}\label{outregdualcond2}
\lim_{t \to \infty}e_{y_i}(t) = \lim_{t \to \infty}(y_i(t) - r(t)) = 0
\end{equation}
for $i = 1,\ldots,N$, all global reference inputs generated by \eqref{smodelg}, all local disturbances represented by \eqref{smodel}, independent from the initial values of the plant \eqref{drs} and of the controller. %Furthermore, the resulting networked controlled MAS has to be stabilized. %This can be seen as a \emph{leader-follower output synchronization problem} in the presence of disturbances.%This output regulation problem can also be seen as a \emph{leader-follower output synchronization problem} in the presence of disturbances, in which the reference model \eqref{smodelg} is the leader and the agents \eqref{drs} are the followers.

%An interesting specialization of the cooperative output regulation problem is  obtained by omitting the global reference model \eqref{smodelg}, i.\:e., by not specifying the reference input. This results in a \emph{leaderless output synchronization problem} in the presence of disturbances, i.\:e., designing a stabilizing cooperative output feedback regulator such that
%\begin{equation}\label{outregdualsync}
%\lim_{t \to \infty}(y_i(t) - y_j(t)) = 0, \quad i,j = 1,\ldots,N
%\end{equation}
%holds for all local disturbances represented by \eqref{smodel} independent from the initial values of the plant \eqref{drs}, of the signal model \eqref{smodel} and of the controller. For this problem, the information exchange through the network is utilized so that the agents can negotiate a common synchronization trajectory for their outputs $y_i$. More specifically, there exists a nontrivial common \emph{synchronization trajectory} $y_s$ such that
%\begin{equation}\label{outregdualsync2}
% \lim_{t \to \infty}(y_i(t) - y_s(t)) = 0, \quad i = 1,\ldots,N.
%\end{equation}
%Different from the cooperative output regulation, the synchronization trajectory $y_s$ is not specified a priori, but depends on the network topology and on the ICs of the agents.

%###########################################################
\section{Design of the cooperative reference observer}\label{sec:netob}
\label{sec:refnw}

%###########################################################
%\subsection{Cooperative reference observer}\label{sec:netob}
Since the global reference model \eqref{smodelg} only transmits the reference input $r$ to the informed agents to keep the communication load small, they estimate the state $w$ of \eqref{smodelg} with the \emph{local reference observers}
\begin{equation}\label{lrefobs}
\dot{\hat{w}}_i(t) = S\hat{w}_i(t) + l_ia_{i0}(r(t) - \hat{r}_i(t)), \quad  t \in \mathbb{R}^+_{t_0},
\end{equation}
in which $\hat{r}_i(t) = p\t\hat{w}_i(t)$, $t \in \overline{\mathbb{R}}^+_{t_0}$, for $i = 1,\ldots,n$, $l_i \in \mathbb{R}^{n_w}$ are the \emph{observer gains} and the observer ICs are $\hat{w}_i(t_0) = \hat{w}_{0,i} \in \mathbb{R}^{n_w}$. Therein, the constants $a_{i0} > 0$, $i = 1,\ldots,n$, describe the broadcast of the output $r$ w.r.t. agent $0$, i.\:e., the global reference model \eqref{smodelg}, to the informed agents. Since $(p\t,S)$ is observable, there always exists an observer gain $l_i$ to ensure an exponentially convergent observer \eqref{lrefobs}, i.\:e., $S-l_ia_{0i}p\t$, $i = 1,\ldots,n$, are Hurwitz matrices. For the uninformed agents, the network has to be utilized to distribute the information about the reference input to the local reference models. Hence, one has to consider the \emph{cooperative reference observer}
\begin{equation}\label{disrobsss}
\dot{\hat{w}}_i(t) = S\hat{w}_i(t) + l_w\textstyle\sum_{j = 1}^Na_{ij}(\hat{r}_j(t) - \hat{r}_i(t))
\end{equation}	 
with $\hat{r}_i(t) = p\t\hat{w}_i(t)$, $t \in \overline{\mathbb{R}}^+_{t_0}$, for $i = n+1,\ldots,N$, to estimate the state $w$ for the uninformed agents. Therein, \eqref{disrobsss} is defined on $t \in \mathbb{R}^+_{t_0}$, $l_w \in \mathbb{R}^{n_w}$ is the \emph{common observer gain} and the observer IC is $\hat{w}_i(t_0) = \hat{w}_{0,i} \in \mathbb{R}^{n_w}$.
%In what follows the observer gain $l_w \in \mathbb{R}^{n_w}$ is determined such that \emph{leader-follower consensus}
%\begin{equation}\label{lfrcons}
% \lim_{t \to \infty}(\hat{r}_i(t) - r(t)) = 0, \quad i = n+1,\ldots,N
%\end{equation}
%is achieved for the uninformed agents. In this problem, the uninformed agents have to follow the output of the reference model \eqref{smodelg}, which is the leader.
%
%\begin{rem}
%Since $(p\t,S)$ is observable by assumption, the property \eqref{lfrcons} directly implies also \emph{state consensus}, i.\:e.,
%\begin{equation}\label{lfrconsst}
%\lim_{t \to \infty}(\hat{w}_i(t) - w(t)) = 0, \quad i = n+1,\ldots,N
%\end{equation}	
%so that the state $w$ of \eqref{smodelg} is estimated by the distributed reference observers \eqref{disrobs}.
%\end{rem}	
By introducing the \emph{observer errors} $e_{w_i} = w - \hat{w}_i$, $i = 1,\ldots,N$, for the local and cooperative reference observer, the \emph{observer error dynamics} related to \eqref{disrobsss} read
\begin{equation}\label{obserrw}
 \dot{e}_{w_i}(t) = Se_{w_i}(t) - \textstyle\sum_{j = 1}^Na_{ij}l_wp\t(\hat{w}_j(t) - \hat{w}_i(t))
\end{equation}
for $i = n+1,\ldots,N$, in view of \eqref{sigmodssg}. In order to simplify \eqref{obserrw}, consider the \emph{Laplacian matrix} $L_{\mathcal{G}}$ associated to the digraph $\mathcal{G}$ describing the communication topology. For the network of the parabolic agents one obtains
\begin{equation}\label{Lmat}
L_{\mathcal{G}} = \begin{bmatrix}
0                 & 0_{1 \times n} & 0_{1 \times N-n}\\
L_{21}            & L_{22}         & 0_{n \times N -n}\\
0_{N-n \times 1}  & L_{32}         & L_{33}	
\end{bmatrix},
\end{equation} 
in which $L_{21} = -\operatorname{col}(a_{10},\ldots,a_{n0})$ and $L_{22} = \operatorname{diag}(a_{10},\linebreak\ldots,a_{n0})$ holds. The first zero row in \eqref{Lmat} results from the fact that the global reference model \eqref{smodelg} (i.\:e., agent $0$) is the root of the communication network. Furthermore, the matrix $L_{21}$ describes the communication of the reference input to the informed agents $i$, $i = 1,\ldots,n$. Hence, no information from the other agents is necessary so that the zero matrix $0_{n \times N -n}$ appears right of $L_{22}$. In contrast, the matrix $0_{N-n \times 1}$ below $L_{21}$ indicates that the uninformed agents $i$, $i = n+1,\ldots,N$, have no access to the reference input. Finally, the matrix $L_{22}$ ensures $L_{\mathcal{G}}1_{N+1} = 0$ with $1_{N+1} = \col{1,\ldots,1} \in \mathbb{R}^{N+1}$, which is a general property of any Laplacian matrix (see, e.\:g., \cite[Lem. 6.2]{Bu19}). With these preparations as well as defining the \emph{local estimate} $\hat{w}^l = \col{\hat{w}_1,\ldots,\hat{w}_n}$ and the \emph{global (cooperative) estimate} $\hat{w}^g = \col{\hat{w}_{n+1},\ldots,\hat{w}_N}$, the aggregated error dynamics \eqref{obserrw} 
%can be rewritten as
%\begin{multline}\label{obserrw2}
% \dot{e}_{w_i}(t) = Se_{w_i}(t) + (e_i\t L_{32} \otimes l_wp\t)\hat{w}^l(t)\\
%  + (e_i\t L_{33} \otimes l_wp\t)\hat{w}^g(t),
%\end{multline}	
%in which $e_i \in \mathbb{R}^{N-n}$ is the $i$-th unit vector. Hence, the aggregated dynamics 
for the observer error $e_w^g = w^g - \hat{w}^g$ with $w^g = (1_{N-n} \otimes I_{n_w})w$ results as
\begin{equation}\label{obserrpre}
 \dot{e}^g_{w}(t) \! = \! F_ge_w^g(t) \! +  \! (L_{32} \otimes l_wp\t)\hat{w}^l(t) \! + \! (L_{33} \otimes l_wp\t)w^g(t),
\end{equation}
in which $F^g = I_{N-n} \otimes S - L_{33} \otimes l_wp\t$ (for the definition and properties of the \emph{Kronecker product} $\otimes$ see, e.\:g., \cite[Ch. 8.2.1]{Bu19}) . This expression can be simplified by observing that \eqref{Lmat} and $L_{\mathcal{G}}1_{N+1} = 0$ imply $L_{33}1_{N-n} = - L_{32}1_n$ so that with $w^l = (1_n \otimes I_{n_w})w$ one obtains $(L_{33} \otimes l_wp\t)w^g =  -(L_{32} \otimes l_wp\t)w^l$.
%\begin{align}\label{nrobserr}
% (L_{33} \otimes l_wp\t)w^g(t) &= (L_{33} \otimes l_wp\t)(1_{N-n} \otimes I_{n_w})w(t)\nonumber\\
%                              &= (L_{33}1_{N-n} \otimes l_wp\t)w(t)\nonumber\\
%                              &= -(L_{32}1_n \otimes l_wp\t)w(t) \nonumber\\
%                              &= -(L_{32} \otimes l_wp\t)w^l(t).
%\end{align}
Hence, inserting this in \eqref{obserrpre} and by defining $e^l_w = w^l - \hat{w}^l$ the error dynamics result as
\begin{subequations}\label{robserragg}
	\begin{align}
	     \dot{e}_w^l(t) &= \operatorname{diag}(S\!-\!l_1a_{10}p\t,\ldots,S\!-\!l_na_{n0}p\t)e_w^l(t)\label{ewsys}\\
	   \dot{e}^g_{w}(t) &= F^ge^g_{w}(t) - (L_{32} \otimes l_wp\t)e_w^l(t).\label{ewgsys}
	\end{align}
\end{subequations}%

The next theorem presents conditions for the stabilization of \eqref{robserragg} and provides a systematic design procedure.
\begin{thm}\emph{\textbf{(Cooperative reference observer)}}\label{theo:robsco}
Assume that the digraph $\mathcal{G}$ associated to the Laplacian matrix $L_{\mathcal{G}}$ in \eqref{Lmat} is connected with the node $0$ as root and let $S-l_ia_{i0}p\t$, $i = 1,\ldots,n$, be Hurwitz matrices. Then, there exists an observer gain $l_w \in \mathbb{R}^{n_w}$ in \eqref{disrobsss} such that $\lim_{t \to \infty}(\hat{w}_i(t) - w(t)) = 0$, $i = 1,\ldots,N$, is achieved for all ICs $w(t_0), \hat{w}_i(t_0) \in \mathbb{R}^{n_w}$, $i = 1,\ldots,N$. A possible choice for the observer gain to ensure this property is $l_w = Qp$ with $Q$ the positive definite solution of the algebraic Riccati equation $SQ + QS\t - 2\nu Qp p\t Q + aI = 0$, where $a > 0$ and $\nu$ such that $\operatorname{Re}\lambda \geq \nu > 0$, $\forall \lambda \in \sigma(L_{33})$.
\end{thm}
The proof of this theorem can be readily obtained from the corresponding results in \cite{Su12} and \cite[Ch. 5.4]{Is17}.

\section{Local state feedback regulator}\label{sec:sfeed}
In this section the \emph{local state feedback regulator}
\begin{align}\label{uxfeed}
 u_i(t) &= \textcolor{blue}{\tfrac{1}{b_i}\big(}-k\t_{w,i}w(t) - k_{v,i}\t(t)v_i(t) - k_{1,i}(t)x_i(1,t)\nonumber\\
        & - \textstyle\int_0^1k_{x,i}(\zeta,t)x_i(\zeta,t)\d \zeta\big)
        = \mathcal{K}_i[w(t),v_i(t),x_i(t),t]
\end{align}
for $i = 1,\ldots,N$ with the \emph{feedback gains} $k_{w,i} \in \mathbb{R}^{n_w}$, $k_{v_i}(t) \in \mathbb{R}^{n_{v_i}}$ and $k_{1,i}(t), k_{x,i}(z,t) \in \mathbb{R}$ is determined, in order to ensure \eqref{outregdualcond2}. The regulator design is based on mapping the agents into backstepping coordinates. This has the advantage that the regulator equations to be solved for determining $k_{w,i}\t$ and $k_{v,i}\t(t)$ in \eqref{uxfeed} can be formulated for the corresponding time-invariant target system. Hence, they are time-invariant, too, facilitating the regulator design. To this end, consider the invertible \emph{backstepping transformation}
\begin{align}\label{btrafo}
 \tilde{x}_i(z,t) &= x_i(z,t) - \textstyle\int_0^zk_i(z,\zeta,t)x_i(\zeta,t)\d\zeta\nonumber\\
                  &= \mathcal{T}_{c,i}[x_i(t),t](z)
\end{align}
for $i = 1,\ldots,N$ with the kernel $k_i(z,\zeta,t) \in \mathbb{R}$. Taking \eqref{smodel} into account, applying \eqref{uxfeed} to \eqref{drs} and introducing the abbreviations $\tilde{g}\t_{k,i} = g\t_{k,i}P_i$, $k = 1,\ldots,4$, results in the closed-loop system
\begin{subequations}\label{clsys}
	\begin{align}
	\dot{w}(t) &= Sw(t)\\
	\dot{v}_i(t) &= S_iv_i(t)\\
	\dot{x}_i(z,t) &= \textcolor{blue}{\lambda_i}x_i''(z,t)\!  + \!a_i(z,t)x_i(z,t)\!  + \! \tilde{g}\t_{1,i}(z,t)v_i(t)\label{xeq2}\\
	x_i'(0,t) &= q_{i}(t)x_i(0,t) + \tilde{g}\t_{2,i}(t)v_i(t)\\
	x_i'(1,t) &= \mathcal{K}_i[w(t),v_i(t),x_i(t),t] +  \tilde{g}\t_{3,i}(t)v_i(t)\label{xbcz0}\\
	          e_{y_i}(t) &= \textcolor{blue}{c_i}x_i(0,t) - p\t w(t) + \tilde{g}\t_{4,i}(t)v_i(t),\label{eyout}
	\end{align}
\end{subequations}
in which $e_{y_i}$ is the \emph{tracking error} defined in \eqref{outregdualcond2}. In order to determine the regulator equations, use \eqref{btrafo}, the feedback gains
\begin{subequations}\label{Kgains}
\begin{align}
\!\!\!\!\!\!k_{1,i}(t) \!\! &=\! -k_i(1,1,t)\textcolor{blue}{-q^{\ell}_i(t)}, k_{x,i}(\zeta,t) \!=\! -k_{i,z}(1,\zeta,t)\\
 k_{w,i}\t &= -\pi_{i,z}\t(1), \quad k_{v,i}\t(t) = \tilde{g}\t_{3,i}(t)- \varphi_{i,z}\t(1,t)
\end{align}	 
\end{subequations}
and the transformation
\begin{equation}\label{tserrdef}
\tilde {\varepsilon}_i(z,t) = \tilde{x}_i(z,t) - \pi_i\t(z)w(t) - \varphi\t_i(z,t)v_i(t)
\end{equation}
for $i = 1,\ldots,N$, with $\pi_i(z) \in \mathbb{R}^{n_w}$ and $\varphi_i(z,t) \in \mathbb{R}^{n_{v_i}}$
%\begin{subequations}\label{Kgains}
%	\begin{align}\label{K1gain}
% 
%	\end{align}
%\end{subequations}
to map \eqref{clsys} into the \emph{target system}
\begin{subequations}\label{clsystt}
	\begin{align}
	\dot{w}(t) &= Sw(t), &&\hspace{-2cm}  t \in \mathbb{R}^+_{t_0}\label{wmodt}\\
	\dot{v}_i(t) &= S_iv_i(t), &&\hspace{-2cm}  t \in \mathbb{R}^+_{t_0}\label{vmodt}\\
	\dot{\tilde{\varepsilon}}_i(z,t) &= \textcolor{blue}{\lambda_i}\tilde{\varepsilon}''_i(z,t) - \mu_i\tilde{\varepsilon}_i(z,t)\label{xeq3tt}\\
	\tilde{\varepsilon}'_i(0,t) &= 0, &&\hspace{-2cm}  t \in \mathbb{R}^+_{t_0}\\
	\tilde{\varepsilon}'_i(1,t) &= 0, &&\hspace{-2cm}  t \in \mathbb{R}^+_{t_0}\label{statetarbc2}\\
	e_{y_i}(t) &=  \textcolor{blue}{c_i}\tilde{\varepsilon}_i(0,t), &&\hspace{-2cm}  t \in  \overline{\mathbb{R}}^+_{t_0}\label{eytt}
	\end{align}
\end{subequations}
for $i = 1,\ldots,N$, and \eqref{xeq3tt} defined on $(z,t) \in (0,1) \times \mathbb{R}^+_{t_0}$. Therein, $\tilde{\varepsilon}_i(z,t)$ is the error between the  transformed state $\tilde{x}_i(z,t)$ and the corresponding reference trajectory $\pi_i\t(z)w(t) + \varphi\t_i(z,t)v_i(t)$ (see \eqref{tserrdef}). The latter is the closed-loop solution, which ensures \eqref{outregdualcond2}. Hence, output regulation is achieved if the tracking error dynamics \eqref{xeq3tt}--\eqref{statetarbc2} is stable (see \eqref{eytt}).
%\begin{subequations}\label{clsyst}
%	\begin{align}
% 	\dot{w}(t) &= Sw(t)\label{wmod}\\
%	\dot{v}_i(t) &= S_iv_i(t)\label{vmod}\\
%	\dot{\tilde{x}}_i(z,t) &= \tilde{x}''_i(z,t) - \mu_i\tilde{x}_i(z,t) + h_{1,i}\t(z,t)v_i(t)\label{xeq3t}\\
%	\tilde{x}'_i(0,t) &= \tilde{g}\t_{2,i}(t)v_i(t)\\
%	\tilde{x}'_i(1,t) &= -k_{w,i}\t w(t) + (\tilde{g}\t_{3,i}(t) -r\t_{v,i}(t))v_i(t)\label{statetarbc}\\
%	e_{y_i}(t) &= \tilde{x}_i(0,t) + \tilde{g}\t_{4,i}(t)v_i(t) - p\t w(t)\label{eyt}
%	\end{align}
%\end{subequations}
Using similar calculations as in \cite{Me13} it can be verified that $k_i(z,\zeta,t) \in \mathbb{R}$ in \eqref{btrafo} has to satisfy the \emph{kernel equations}
\begin{subequations}\label{keq}
	\begin{align}
	 k_{i,t}(z,\zeta,t) &= \textcolor{blue}{\lambda_i}k_{i,zz}(z,\zeta,t) - \textcolor{blue}{\lambda_i}k_{i,\zeta\zeta}(z,\zeta,t)\nonumber\\
	                & \qquad - (a_i(\zeta,t) + \mu_i)k_i(z,\zeta,t)\label{kpde}\\
	 k_{i,\zeta}(z,0,t) &= q_i(t)k_i(z,0,t)\label{kBC1}\\
  	 k_i(z,z,t) &= q_i(t) - \textstyle\int_0^z\tfrac{a_i(\zeta,t) + \mu_i}{2\textcolor{blue}{\lambda_i}}\d\zeta\label{kBC3}.	
	\end{align}
\end{subequations}
Therein, $\mu_i \in \mathbb{R}$ and \eqref{kpde} is defined on $(z,\zeta,t) = \{0 < \zeta < z < 1\} \times \mathbb{R}_{t_0}^+$. It is shown in \cite{Vaz08,Me09,Me13} that the kernel equations \eqref{keq} have a strong solution, which is in the Gevrey class $G_{\alpha}(\mathbb{R}^+_{t_0})$, $1 \leq \alpha < 2$, w.r.t. time. The inverse transformation $\mathcal{T}^{-1}_{c,i}$ is also a backstepping transformation with a kernel following from similar kernel equations (see \cite{Vaz08,Me09,Me13}).
%In order to trace the output regulation problem back to a stabilization problem, the signal models \eqref{wmod} and \eqref{vmod} are decoupled from the controlled agent dynamics \eqref{xeq3t}--\eqref{statetarbc}. This is achieved with the
%%\begin{subequations}\label{orpgains}
%%	\begin{align}
%%	\end{align}
%%\end{subequations}
%the corresponding target system is given by
Differentiating \eqref{tserrdef} w.r.t. time, inserting the system \eqref{clsys} mapped in the backstepping coordinates \eqref{btrafo} yields  \eqref{clsystt} if $\pi_i\t(z)$, $i = 1,\ldots,N$, is the solution of the \emph{reference regulator equations}
\begin{subequations}\label{greq}
	\begin{align}
     \!\!\!\!\textcolor{blue}{\lambda_i}\pi_{i,zz}\t(z) \!-\! \mu_i\pi_i\t(z) \!-\! \pi_i\t(z) S &= 0\t,\quad  z \in (0,1)\label{greqode}\\
                                   \pi_{i,z}\t(0) &= 0\t\label{greqic1}\\
                                           \textcolor{blue}{c_i}\pi_{i}\t(0) &= p\t\label{greqic2}
	\end{align}
\end{subequations}
and the $\varphi_i\t(z,t)$, $i = 1,\ldots,N$, is the solution of the \emph{disturbance regulator equations}
\begin{subequations}\label{lreq}
	\begin{align}
	\varphi_{i,t}\t(z,t) &= \textcolor{blue}{\lambda_i}\varphi_{i,zz}\t(z,t) - \mu_i\varphi_i\t(z,t)- \varphi_i\t(z,t) S\nonumber\\
	 & \quad + h_{1,i}\t(z,t),\quad  (z,t) = (0,1) \times \mathbb{R}^+_{t_0}\label{greqpde}\\
	\varphi_{i,z}\t(0,t) &= \tilde{g}_{2,i}\t(t), &&\hspace{-3cm} t \in \mathbb{R}^+_{t_0} \label{lreqbc1}\\
	 \textcolor{blue}{c_i}\varphi_{i}\t(0,t) &= -\tilde{g}_{4,i}\t(t), &&\hspace{-3cm} t \in \mathbb{R}^+_{t_0}\label{lreqbc2}
	\end{align}
\end{subequations}
with $h_{1,i}\t(z,t) = k_i(z,0,t)\textcolor{blue}{\lambda_i}\tilde{g}_{2,i}\t(t) + \mathcal{T}_{c,i}[\tilde{g}\t_{1,i}(t),t](z)$. The next lemma shows that the regulator equations \eqref{greq} admit an explicit solution.
\begin{lem}\emph{\textbf{(Reference regulator equations)}}\label{lem:lreq}
The regulator equations \eqref{greq} have the strong solution
\begin{equation}
 \!\pi_i(z) \! = \! \begin{bmatrix}
 I_{n_w} & 0
 \end{bmatrix}\! \exp\left\{\!\begin{bmatrix}
 0 & I_{n_w}\\
 \textcolor{blue}{\frac{1}{\lambda_i}}(\mu_iI + S\t) & 0
 \end{bmatrix}\!\!z\!\right\}\!\!\begin{bmatrix}\textcolor{blue}{\tfrac{1}{c_i}}p \\ 0\end{bmatrix}
\end{equation}
for $i = 1,\ldots,N$.
\end{lem}
This result directly follows from the reformulation of \eqref{greq} as an initial value problem for $\pi_i(z)$ and $\pi_{i,z}(z)$ and from the corresponding solution.
%\begin{pf}
%By introducing $\phi_i = \col{\pi_i,\pi_{i,z}}$, $i = 1,\ldots,N$, the ODE \eqref{greqode} can be represented by 
%$\phi_{i,z}(z) = A_i\phi_i(z)$, $z \in (0,1]$. With \eqref{greqic1} and \eqref{greqic2} this yields an initial value problem (IVP) so that from its solution the strong solution can easily be deduced.	 \hfill $\Box$
%\end{pf}
Much more challenging is to solve the disturbance regulator equations \eqref{lreq}, because this is a \emph{Cauchy problem} for the PDE \eqref{greqpde}. For this, the next lemma presents a systematic solution procedure establishing the solvability of \eqref{lreq}.
\begin{lem}\emph{\textbf{(Disturbance regulator equations)}}\label{lem:lreq}
A strong solution of the regulator equations \eqref{lreq} is
 \begin{multline}\label{lreqsol}
	\varphi_i(z,t) = 
	\begin{bmatrix}
	I_{n_{v_i}} & 0
	\end{bmatrix}
	\Big(\sum_{j=0}^{\infty}\Phi_{j,i}(z)\d_t^j
	\begin{bmatrix}
	- \textcolor{blue}{\tfrac{1}{c_i}}\tilde{g}_{4,i}(t)\\
	\tilde{g}_{2,i}(t)
	\end{bmatrix}\\
	 -\tfrac{1}{\textcolor{blue}{\lambda_i}}\int_0^z\sum_{j=0}^{\infty}\Phi_{j,i}(\zeta)
	\begin{bmatrix}
	0 & I_{n_{v_i}}
	\end{bmatrix}\t\d_t^jh_{1,i}(z-\zeta,t)\d\zeta\Big)
\end{multline}
for $i = 1,\ldots,N$ with $\Phi_{j,i}(z) = \begin{bmatrix}
 \phi^i_{1,j}(z) & \ldots & \phi^i_{2n_{v_i},j}(z)
 \end{bmatrix}$, $j \geq 0$, in which the elements $\phi^i_{k,j}(z)$, $k = 1,\ldots,2n_{v_i}$, result from the recursion $\phi^i_{k,0}(z) = \e^{A_{0,i}z}e_k$, 
$\phi^i_{k,j}(z) = \tint_0^z\e^{A_{0,i}(z-\zeta)}A_{1,i}\phi^i_{k,j-1}(\zeta)\d\zeta$, $j \geq 1$ and
\begin{equation}
 A_{0,i} = \begin{bmatrix}
  0 & I_{n_{v_i}}\\
   \textcolor{blue}{\frac{1}{\lambda_i}}(\mu_iI + S_i\t)\;\; & 0
 \end{bmatrix}, \;
 A_{1,i} = \begin{bmatrix}
 0 & 0_{n_{v_i}}\\
 \textcolor{blue}{\frac{1}{\lambda_i}}I_{n_{v_i}} & 0
 \end{bmatrix}.
\end{equation}
\end{lem}
\begin{pf}
Introduce $\phi_i = \col{\varphi_i,\varphi_{i,z}}$, $i = 1,\ldots,N$, and formally apply the Laplace transform to \eqref{lreq} so that the PDE \eqref{greqpde} can be represented by 
$\check{\phi}_{i,z}(z,s) = A_i(s)\check{\phi}_i(z,s) + [0 \;\; I_{n_w}]\t \check{h}_{1,i}(z,s)$, $z \in (0,1]$, with $A_i(s) = A_{1,i}s + A_{0,i}$. For a mathematical justification of this formal approach in the sense of the Mikusinski's operator calculus, the reader is referred to \cite{Ru08}. With \eqref{lreqbc1} and \eqref{lreqbc2} this yields an IVP, which can be solved by making use of the matrix exponential $\e^{A_i(s)z}$. In order to determine the latter, consider the solution $\psi_i(z,s) = \e^{A_i(s)z}\psi_{i,0}(s)$ of the IVP $\psi_{i,z}(z,s) = A_i(s)\psi_i(z,s), \quad z \in (0,1], \psi_i(0,s) = \psi_{i,0}(s)$.
Let $\psi_i(z,s) = \sum_{j=0}^{\infty}\Phi_{j,i}(z)s^j$ and insert this in the IVP so that equating coefficients w.r.t. $s$ directly yields $\d_z\Phi_{0,i}(z) = A_{0,i}\Phi_{0,i}(z)$,
$\d_z\Phi_{j,i}(z) = A_{0,i}\Phi_{j,i}(z) + A_{1,i}\Phi_{j-1,i}(z)$, $j \geq 1$, for the coefficients in $\psi_i(z,s)$. In order to determine the $k$-th column of the matrix exponential $\e^{A_i(s)z}$, introduce the $k$-th unit vector $e_k \in \mathbb{R}^{2n_{v_i}}$ and let $\psi_{0,i}(s) = e_k$. Then, the series representation for $\psi_i(z,s)$ implies $\Phi_{0,i}(0) = e_k$ and $\Phi_{j,i}(0) = 0$, $j \geq 1$. This and the recursion for $\Phi_{j,i}(z)$ determines the recursion in the theorem to compute the $k$-th column of $\e^{A_i(s)z}$. By applying Theorem 1 in \cite{Ru08} to $A_i(s)$ in the IVP it is straightforward to verify that the series in \eqref{lreqsol} converge uniformly on $(z,t) \in [0,1] \times  \mathbb{R}^+_{t_0}$ if the elements of $\tilde{g}_{j,i}(t)$, $j = 2,4$, and of $h_{1,i}(z,t)$ are in $G_{\alpha}(\mathbb{R}^+_{t_0})$, $1 \leq \alpha < 2$, w.r.t. time. The former conditions are fulfilled by assumption (see Section \ref{sec:probform}). Since, in addition, the kernel $k_i(z,\zeta,t)$ is in  $G_{\alpha}(\mathbb{R}^+_{t_0})$, $1 \leq \alpha < 2$, w.r.t. time, the expression for $h_{1,i}\t(z,t)$ implies the same for $h_{1,i}(z,t)$. This shows that \eqref{lreqsol} is the strong solution of \eqref{lreq}. \hfill $\Box$
\end{pf}

The next theorem presents the stability result for the tracking error dynamics, which implies output regulation.
\begin{thm}\emph{\textbf{(Local state feedback regulator)}}\label{thm:tstab}
Assume that $\mu_i > 0$, $i = 1,\ldots,N$. Let $k_i(z,\zeta,t)$, $\pi_i\t(z)$  and $\varphi_i\t(z,t)$ be the solutions of the kernel equations \eqref{keq} and of the regulator equations \eqref{greq} and \eqref{lreq}. Then, the state feedback regulator \eqref{uxfeed} with the feedback gains \eqref{Kgains} achieves output regulation \eqref{outregdualcond2}. The dynamics of the tracking error $e_i(t) = \{e_i(z,t), z \in [0,1]\}$ with $e_i(z,t) = x_i(z,t) - \mathcal{T}_{c,i}^{-1}[\pi_i\t,t](z)w(t) - \mathcal{T}_{c,i}^{-1}[\varphi_i\t(\cdot,t),t](z)v_i(t)$ are uniformly exponentially stable in the $L_2$-norm, i.\:e., $\|e_i(t)\|_{L_2} \leq M_i\operatorname{e}^{-\mu_i (t-t_0)}\|e_i(t_0)\|_{L_2},\quad t \in \overline{\mathbb{R}}^+_{t_0}$ 
for all $e_i(t_0) \in H^2(0,1)$, $i = 1,\ldots,N$, satisfying the BCs of the tracking error dynamics, an $M_i \geq 1$ and any $t_0 \geq 0$.
\end{thm}	
The tracking error dynamics \eqref{xeq3tt}--\eqref{statetarbc2} coincides with the corresponding dynamics in \cite{Deu15a} with the BC at $z = 0$ being of Neumann type. Hence, the proof of Theorem \ref{thm:tstab} can be readily deduced from the corresponding result in \cite{Deu15a}.

\section{Local disturbance observer}\label{sec:dlocobs}
In order to estimate the states of the local disturbance model \eqref{smodel} and the individual agents \eqref{drs}, the \emph{local disturbance observer}
\begin{subequations}\label{ldobs}
\begin{align}
   \dot{\hat{v}}_i(t) &= S_i\hat{v}_i(t) + l_{v_i}(t)(\eta_i(t) - \textcolor{blue}{c_{m,i}}\hat{x}_i(1,t)) &&\label{dobs1}\\
 \dot{\hat{x}}_i(z,t) &= \textcolor{blue}{\lambda_i}\hat{x}''_i(z,t) + a_i(z,t)\hat{x}_i(z,t) +  \tilde{g}_{1,i}\t(z,t)\hat{v}_i(t)\nonumber\\
                      &\quad  + l_{x_i}(z,t)(\eta_i(t) - \textcolor{blue}{c_{m,i}}\hat{x}_i(1,t))\label{dobsss}\\
      \hat{x}'_i(0,t) &= q_i(t)\hat{x}_i(0,t) + \tilde{g}_{2,i}\t(t)\hat{v}_i(t), &&\hspace{-2cm} t \in \mathbb{R}^+_{t_0}\label{obsrb2}\\
   \hat{x}'_i(1,t) &=  \textcolor{blue}{\tfrac{q^{\ell}_i(t)}{c_{m,i}}\eta_i(t)} +  \textcolor{blue}{b_i}u_i(t) + \tilde{g}_{3,i}\t(t)\hat{v}_i(t)\nonumber\\
                      & \quad + l_{1,i}(t)(\eta_i(t) - \textcolor{blue}{c_{m,i}}\hat{x}_i(1,t)),&& \hspace{-2cm} t \in \mathbb{R}^+_{t_0}\label{dobsu}  
	 \end{align}
\end{subequations}
is designed with \eqref{dobs1} defined on $t \in \mathbb{R}^+_{t_0}$ and \eqref{dobsss} on  $(z,t) \in (0,1) \times \mathbb{R}^+_{t_0}$ as well as the ICs $\hat{v}_i(t_0) = \hat{v}_{i,0} \in \mathbb{R}^{n_{v_i}}$ and $\hat{x}_i(z,t_0) = \hat{x}_{i,0}(z) \in \mathbb{R}$, $z \in [0,1]$. Note that \eqref{obsrb2} does not need an output injection term in the collocated setting for mapping the corresponding  observer error dynamics \eqref{ldobs} into a stable target system. This was already shown for the time-invariant case in \cite{Sm05a}. In the sequel, the \emph{observer gains} $l_{v_i}(t) \in \mathbb{R}^{n_{v_i}}$, $l_{x_i}(z,t) \in \mathbb{R}$ and $l_{1,i}(t) \in \mathbb{R}$ are determined to ensure that the corresponding \emph{observer error dynamics} 
\begin{subequations}\label{obserrdyn}
\begin{align}
   \dot{e}_{v_i}(t) &= S_ie_{v_i}(t) - l_{v_i}(t)\textcolor{blue}{c_{m,i}}e_{x_i}(1,t)\label{everr}\\
 \dot{e}_{x_i}(z,t) &= \textcolor{blue}{\lambda_i}e''_{x_i}(z,t) + a_i(z,t)e_{x_i}(z,t) + \tilde{g}_{1,i}\t(z,t)e_{v_i}(t)\nonumber\\
                    & \quad -l_{x_i}(z,t)\textcolor{blue}{c_{m,i}}e_{x_i}(1,t)\label{errs2}\\
  e'_{x_i}(0,t) &= q_i(t)e_{x_i}(0,t) + \tilde{g}_{2,i}\t(t)e_{v_i}(t)\label{errrb1}\\	                       
           e'_{x_i}(1,t) &= \tilde{g}_{3,i}\t(t)e_{v_i}(t) - l_{1,i}(t)\textcolor{blue}{c_{m,i}}e_{x_i}(1,t)\label{errs4}	 
\end{align}
\end{subequations}
with $e_{v_i}(t) = v_i(t) - \hat{v}_i(t)$ and $e_{x_i}(z,t) = x_i(z,t) - \hat{x}_i(z,t)$ are stabilized. The observer error dynamics \eqref{obserrdyn} are mapped into the \emph{PDE-ODE cascade}
\begin{subequations}\label{tfirst2}
	\begin{align}
	\dot{\varepsilon}_{x_i}(z,t) &= \textcolor{blue}{\lambda_i}\varepsilon''_{x_i}(z,t) - \bar{\mu}_i\varepsilon_{x_i}(z,t)\label{errs2t2}\\
	\varepsilon'_{x_i}(0,t) &=0\label{errrb1t2}\\	            
	\varepsilon'_{x_i}(1,t) &= 0\label{errs4t2}\\
	\dot{e}_{v_i}(t) &= (S_i - l_{v_i}(t)\textcolor{blue}{c_{m,i}}\gamma_i\t(1,t))e_{v_i}(t)\nonumber\\ 
	                 & \quad - l_{v_i}(t)\textcolor{blue}{c_{m,i}}\varepsilon_{x_i}(1,t)\label{everrt2}		 
	\end{align}
\end{subequations}
%\begin{subequations}\label{tfirst}
%	\begin{align}
%	\dot{e}_{v_i}(t) &= S_ie_{v_i}(t) - l_{v_i}(t)\tilde{e}_{x_i}(1,t)\label{everrt}\\
%	\dot{\tilde{e}}_{x_i}(z,t) &= \tilde{e}''_{x_i}(z,t) - \bar{\mu}_i\tilde{e}_{x_i}(z,t) + \bar{h}\t_{1,i}(z,t)e_{v_i}(t)\nonumber\\
%	& \quad -\tilde{l}_{x_i}(z,t)\tilde{e}_{x_i}(1,t)\label{errs2t}\\
%	\tilde{e}'_{x_i}(0,t) &=\tilde{g}_{2,i}\t(t)e_{v_i}(t)\label{errrb1t}\\	                       
%	\tilde{e}'_{x_i}(1,t) &= \tilde{g}_{3,i}\t(t)e_{v_i}(t)\label{errs4t}	 
%	\end{align}
%\end{subequations}
by making use of the invertible \emph{backstepping transformation}
\begin{subequations}
\begin{align}
e_{x_i}(z,t) &= \tilde{e}_{x_i}(z,t) - \textstyle\int_z^1p_{i}(z,\zeta,t)\tilde{e}_{x_i}(\zeta,t)\d\zeta\nonumber\\
&= \mathcal{T}^{-1}_{o,i}[\tilde{e}_{x_i}(t),t](z)\label{vtrafo2}\\
\varepsilon_{x_i}(z,t) &= \tilde{e}_{x_i}(z,t) - \gamma_i\t(z,t)e_{v_i}(t)\label{decoupltraf}       
\end{align}
\end{subequations}
with the kernel $p_{i}(z,\zeta,t) \in \mathbb{R}$ and $\gamma_i(z,t) \in \mathbb{R}^{n_{v_i}}$. Therein, the \emph{observer gains}
\begin{subequations}\label{obsgain1}
\begin{align}
 l_{1,i}(t) &= -\textcolor{blue}{\tfrac{1}{c_{m,i}}}p_i(1,1,t)\\
 \!\!\!\!l_{x_i}(z,t) &= \mathcal{T}^{-1}_{o,i}[\gamma_i\t(t)l_{v_i}(t),t](z) - p_{i,\zeta}(z,1,t)\textcolor{blue}{\tfrac{\lambda_i}{c_{m,i}}}	 	 
\end{align}
\end{subequations}
are utilized. From the related derivation in \cite[Ch 8.3]{Me13} it is readily deduced that $p_{i}(z,\zeta,t)$ has to be the solution of the \emph{kernel equations}
\begin{subequations}\label{obsker}
 \begin{align}
  p_{i,t}(z,\zeta,t) &= \textcolor{blue}{\lambda_i}p_{i,zz}(z,\zeta,t) - \textcolor{blue}{\lambda_i}p_{i,\zeta\zeta}(z,\zeta,t)\nonumber\\
                     & \quad + (a_i(z,t) + \bar{\mu}_i)p_i(z,\zeta,t)\label{ckpde2}\\
  p_{i,z}(0,\zeta,t) &= q_i(t)p_i(0,\zeta,t)\\
          p_i(z,z,t) &= q_i(t) - \tint_0^z\tfrac{a_i(\zeta,t) + \bar{\mu}_i}{2\textcolor{blue}{\lambda_i}}\d\zeta\label{ckpde32}
 \end{align}
\end{subequations}
with $\bar{\mu}_i \in \mathbb{R}$ and \eqref{ckpde2} defined on $(z,\zeta,t) = \{0 < z < \zeta < 1\} \times \mathbb{R}_{t_0}^+$. It is shown in \cite[Ch 8.3]{Me13} that they admit a strong solution, which is in the Gevrey class $G_{\alpha}(\mathbb{R}^+_{t_0})$, $1 \leq \alpha < 2$, w.r.t. time. The corresponding inverse transformation $\mathcal{T}_{o,i}$ is also a backstepping transformation, which can be determined by solving kernel equations of the same type as \eqref{obsker} (see \cite[Ch 8.3]{Me13}).
%In order to map \eqref{tfirst} into the \emph{PDE-ODE cascade}
%\begin{subequations}\label{tfirst2}
%	\begin{align}
%	\dot{\varepsilon}_{x_i}(z,t) &= \varepsilon''_{x_i}(z,t) - \bar{\mu}_i\varepsilon_{x_i}(z,t)\label{errs2t2}\\
%           \varepsilon'_{x_i}(0,t) &=0\label{errrb1t2}\\	            
%	\varepsilon'_{x_i}(1,t) &= 0\label{errs4t2}\\
%	\dot{e}_{v_i}(t) &= (S_i - l_{v_i}(t)\gamma_i\t(1,t))e_{v_i}(t) 
%	- l_{v_i}(t)\varepsilon_{x_i}(1,t),\label{everrt2}		 
%	\end{align}
%\end{subequations}
%and the \emph{observer gain}
%\begin{equation}\label{obsgain3}
% \tilde{l}_{x_i}(z,t) = \gamma_i\t(z,t)l_{v_i}(t)
%\end{equation}
%are utilized. 
Differentiating \eqref{decoupltraf} w.r.t. time and inserting system \eqref{obserrdyn} mapped in the backstepping coordinates \eqref{vtrafo2} yields \eqref{tfirst2} if $\gamma_i\t(z,t)$ solves the initial boundary value problem (IBVP) 
\begin{subequations}\label{decoupleq}
	\begin{align}
	 \gamma\t_{i,t}(z,t) &= \textcolor{blue}{\lambda_i}\gamma_{i,zz}\t(z,t) - \bar{\mu}_i\gamma_i\t(z,t)- \gamma_i\t(z,t)S_i \nonumber\\
	 & \quad + \bar{h}_{1,i}\t(z,t), \quad  (z,t) \in (0,1) \times \mathbb{R}^+_{t_0}\label{decplpde}\\
	 \gamma_{i,z}\t(0,t) &= \tilde{g}_{2,i}\t(t), &&\hspace{-4cm}  t \in \mathbb{R}^+_{t_0}\label{debc1}\\
	 \gamma_{i,z}\t(1,t) &= \tilde{g}_{3,i}\t(t), &&\hspace{-4cm}  t \in \mathbb{R}^+_{t_0}\label{debc2}
	\end{align}
\end{subequations}
with constant coefficients facilitating the well-posedness proof. Therein,  $\bar{h}\t_i(z,t) = \mathcal{T}_{o,i}[\tilde{g}_{1,i}\t(\cdot,t),t](z)\linebreak + \mathcal{T}_{o,i}[p_i(\cdot,1,t,),t](z)\textcolor{blue}{\lambda_i}\tilde{g}\t_{3,i}(t)$ holds. The next lemma asserts the solvability of \eqref{decoupleq}.

\begin{lem}\emph{\textbf{(Solvability of \pmb{\eqref{decoupleq}})}}\label{lem:decoupleq}
The IBVP \eqref{decoupleq} for $i = 1,\ldots,N$ have a unique solution in $L_2(0,1)$ for any IC $\gamma_i\t(z,t_0) = \gamma\t_{i,0}(z)$ with the elements of $\gamma_{i,0}$ in $H^2(0,1)$ satisfying the BCs \eqref{debc1} and \eqref{debc2}. Furthermore, the elements of $\gamma_i\t(1,\cdot)$, $i = 1,\ldots,N$, are in $G_{\alpha}(\mathbb{R}^+_{t_0})$, $1 \leq \alpha < 2$.
\end{lem}
For the proof see Appendix \ref{app}. In order to ensure a stable PDE-ODE cascade \eqref{tfirst2}, the $e_{v_i}$-subsystem has to be stabilized by a  suitable observer gain $l_{v,i}(t)$. The next lemma presents the condition for its existence.
\begin{lem}\emph{\textbf{(Observability)}}\label{lem:declode}
Consider $c(t) \in \mathbb{R}^{n_{v_i}}$ and $A(t) \in \mathbb{R}^{n_{v_i} \times n_{v_i}}$ and let their elements be sufficiently smooth. Define the operator $\mathcal{M}_{A}c\t(t) = c\t(t)A(t) + \dot{c}\t(t)$ recursively by
$\mathcal{M}^k_{A}c\t(t) = \mathcal{M}_{A}(\mathcal{M}^{k-1}_{A}c\t(t))$, $k \geq 1$, and  $\mathcal{M}^0_{A}c\t(t) = c\t(t)$. If $\det Q_{o,i}(t) \neq 0$, $t \in \overline{\mathbb{R}}^+_{t_0}$, with $Q_{o,i}(t) = [\mathcal{M}^0_{S_i}\gamma_i(1,t) \;\; \ldots \;\; \mathcal{M}^{n_{v_i}-1}_{S_i}\gamma_i(1,t)]\t$, then there exists a bounded observer gain $l_{v_i}(t)$, $i = 1,\ldots,N$, such that the unforced dynamics of the $e_{v_i}$-system \eqref{everrt2} are uniformly exponentially stable. More precisely, there exists a fundamental matrix $\Psi_i(t,\tau) : \overline{\mathbb{R}}^+_{t_0} \times \overline{\mathbb{R}}^+_{t_0} \to \mathbb{R}^{n_{v_i} \times n_{v_i}}$ being the solution of the IVP $\Psi_{i,t}(t,\tau) = \big(S_i - l_{v_i}(t)\textcolor{blue}{c_{m,i}}\gamma_i\t(1,t)\big)\Psi_i(t,\tau)$, $\Psi_i(\tau,\tau) = I$, and a positive constant $\kappa_i$ such that $\|\Psi_i(t,\tau)\| \leq \kappa_i\e^{-\bar{\mu}_{v_i}(t - \tau)}$ for all $t$ and $\tau$ with $t \geq \tau$.
\end{lem}
\begin{pf}
Since Lemma \ref{lem:decoupleq} ensures that $\gamma_i\t(1,\cdot) \in G_{\alpha}(\mathbb{R}^+_{t_0})$, $1 \leq \alpha < 2$, the matrix $Q_{o,i}(t)$ exists. The condition $\det Q_{o,i}(t) \neq 0$ implies that the pair $(\gamma_i\t(1,t),S_i)$, $i = 1,\ldots,N$, is uniformly observable for $t \in \overline{\mathbb{R}}^+_{t_0}$ (see, e.\:g., \cite[Th. 9.10]{Rug96}). This ensures the existence of an observer gain $l_{v_i}(t)$ such that the unforced dynamics of the $e_{v_i}$-system is uniformly exponentially stable (see, e.\:g., \cite{Ik75}). \hfill $\Box$
\end{pf}%
%\begin{rem}
If the condition of Lemma \ref{lem:declode} is satisfied, then various methods to determine the observer gain $l_{v_i}(t)$ become available. For example, an eigenvalue assignment using the time-variant observer canonical form can be utilized (see, e.\:g., \cite{Val95}). 
%It should be noted that this requires $\gamma_i\t(1,t)$ in explicit form, which is, in general, hard to obtain. Consequently, a numerical approach is presented in Section \ref{sec:ex}, in order to provide a systematic design procedure.
%\end{rem}	
The next theorem asserts the uniform exponential stability of the observer error dynamics \eqref{obserrdyn}.
\begin{thm}\emph{\textbf{(Local disturbance observer)}}\label{th:dobs}
Consider the observer \eqref{ldobs} and let the observer gains  $l_{v_i}(t)$, $l_{x_i}(z,t)$ and $l_{1,i}(t)$, $i = 1,\ldots,N$, be given by \eqref{obsgain1}. Assume that $\bar{\mu}_i > 0$, $\bar{\mu}_{v_i} > 0$ (see Lemma \ref{lem:declode}) and $\bar{\mu}_i \neq \bar{\mu}_{v_i}$ implying $\alpha_{o,i} = \min(\bar{\mu}_i,\bar{\mu}_{v_i}) > 0$. Then, the dynamics of the  observer error $e_{o,i}(t) = \operatorname{col}(e_{v_i}(t),e_{x_i}(t))$ with $e_{x_i}(t) = \{e_{x_i}(z,t), z \in [0,1]\}$ are uniformly exponentially stable in the norm $\|\cdot\|_{X_{o,i}} = (\|\cdot\|^2_{\mathbb{C}^{n_{v_i}}} + \|\cdot\|_{L_2}^2)^{\frac{1}{2}}$, i.\:e., 
$\|e_{o}(t)\|_{X_{o,i}} \leq \bar{M}_{i}\e^{-\alpha_{o,i}(t-t_0)}	
\|e_{o}(t_0)\|_{X_{o,i}}$, $t \in \overline{\mathbb{R}}^+_{t_0}$ for all $e_{o,i}(t_0) \in \mathbb{C}^{n_{v_i}} \oplus H^2(0,1)$ satisfying the BCs of the observer error dynamics, an $\bar{M}_i \geq 1$ and any $t_0 \geq 0$. 
\end{thm}
For the proof see Appendix \ref{app}.

%\subsection{Output feedback regulator}
%The \emph{output feedback regulator} results from  inserting the state estimates in \eqref{uxfeed} yielding
%\begin{equation}\label{outfeedregobs}
% u_i(t) = \mathcal{K}_i[\hat{w}_i(t),\hat{v}_i(t),\hat{x}_i(t),t].
%\end{equation}
%In particular, the estimate $\hat{w}_i$ is obtained from the cooperative reference observer \eqref{lrefobs}, \eqref{disrobs}, the estimates $\hat{x}_i$ and $\hat{v}_i$ of the local agent and disturbance model states are obtained from the disturbance observer \eqref{ldobs}. 

\section{Cooperative output feedback regulator}\label{sec:coufeedreg}
%In what follows, the stability is investigated for the closed-loop networked MAS resulting from applying the cooperative output feedback regulator designed in Sections \ref{sec:refnw} and \ref{sec:locreg} for both setups, which implies cooperative output regulation and leaderless output synchronization, respectively. 
The next theorem shows that the compensator resulting from applying the estimates of the local reference observers \eqref{lrefobs}, the cooperative reference observer \eqref{disrobsss} and the local disturbance observer \eqref{ldobs} in \eqref{uxfeed} ensures cooperative output regulation.
%clarifies the closed-loop stability ensuring cooperative output regulation according to \eqref{outregdualcond2}. In particular, 
%For this, the dynamics of the closed-loop system is represented by the state $x_{\text{co}} = \col{e^l_w,e^g_w,e_v,e_x,\hat{x}-\bar{\pi}_i\t\hat{w}_i 
%-\bar{\varphi}_i\t\hat{v}_i}$, in which $\bar{\pi}_i\t(z,t) = \mathcal{T}_{c,i}^{-1}[\pi\t,t](z)$ and $\bar{\varphi}_i\t(z,t) = \mathcal{T}_{c,i}^{-1}[\varphi_i\t(t),t](z)$ are utilized (see \eqref{btrafo}). Therein, $e^l_w(t) \in \mathbb{R}^{nn_w}$ and $e^g_w(t) \in \mathbb{R}^{(N-n)n_w}$ are the observer errors for the local and the cooperative reference observer (see Section \ref{sec:netob}), $e_v = \col{e_{v_1},\ldots,e_{v_N}}$ and $e_x = \col{e_{x_1},\ldots,e_{x_N}}$ are the observer errors for the parabolic agent disturbance observer and $\hat{x} = \col{\hat{x}_1,\ldots,\hat{x}_N}$ is its state w.r.t. the agent dynamics (see Section \ref{sec:dlocobs}). 
\begin{thm}\emph{\textbf{(Cooperative output regulation)}}\label{thm:sepprinc}
Consider the observers \eqref{lrefobs}, \eqref{disrobsss} and \eqref{ldobs} designed according to Theorems \ref{theo:robsco}, \ref{thm:tstab} and use their estimates in $u_i(t) = \mathcal{K}_i[\hat{w}_i(t),\hat{v}_i(t),\hat{x}_i(t),t]$, $i = 1,\ldots,N$ (see \eqref{uxfeed}). Then, cooperative output regulation \eqref{outregdualcond2} is achieved and the dynamics of the closed-loop state $x_{\text{\emph{co}}} = \col{e^l_w,e^g_w,e_v,e_x,\hat{e}_{x_1},\ldots,\hat{e}_{x_N}}$, in which $\hat{e}_{x_i} = \hat{x}_i - \bar{\pi}_i\t\hat{w}_i -\bar{\varphi}_i\t\hat{v}_i$ with $\bar{\pi}_i\t(z,t) = \mathcal{T}_{c,i}^{-1}[\pi\t,t](z)$ and $\bar{\varphi}_i\t(z,t) = \mathcal{T}_{c,i}^{-1}[\varphi_i\t(t),t](z)$ are utilized (see \eqref{btrafo}) are uniformly exponentially stable in the norm $\|\cdot\|_{X_{\text{\emph{co}}}} = (\|\cdot\|^2_{\mathbb{C}^{nn_w}} + \|\cdot\|^2_{\mathbb{C}^{(N-n)n_w}} + \|\cdot\|^2_{\mathbb{C}^{n_{v_1} + \ldots + n_{v_N}}}  + \underbrace{\|\cdot\|^2_{L_2} + \ldots  +\|\cdot\|^2_{L_2}}_{\text{2N-times}})^{\frac{1}{2}}$. 
\end{thm}
For the proof see Appendix \ref{app}.

%In the case of leaderless output synchronization, the closed-loop system is represented by the state $x_{\text{lls}} = \col{e_{\hat{w}^g},e_v,e_x,\hat{x}-\bar{\pi}_i\t\hat{w}_i 
%-\bar{\varphi}_i\t\hat{v}_i}$, in which different from the cooperative output regulation setup the synchronization error $e_{\hat{w}^g}(t) \in \mathbb{R}^{(N-1)n_w}$ appears instead of the observer errors $e^l_w$ and $e^g_w$ (see Section \ref{sec:synobs}). With this, the next theorem asserts output synchronization.
%\begin{thm}\emph{\textbf{(Leaderless output synchronization)}}\label{thm:outsync}
%Consider the cooperative output feedback regulator \eqref{disrobss}, \eqref{ldobs} and \eqref{outfeedregobs} designed according to Theorems \ref{theo:robscos}, \ref{thm:tstab} and \ref{th:dobs}. Then, leaderless output synchronization \eqref{outregdualsync} is achieved and the dynamics of the closed-loop state $x_{\emph{\text{lls}}}(t)$ are uniformly exponentially stable in the norm $\|\cdot\|_{X_{\emph{\text{lls}}}} = (\|\cdot\|^2_{\mathbb{C}^{(N-1)n_w}} + \|\cdot\|^2_{\mathbb{C}^{n_{v_1} + \ldots + n_{v_N}}} + \underbrace{\|\cdot\|^2_{L_2} + \ldots  +\|\cdot\|^2_{L_2}}_{\text{2N-times}})^{\frac{1}{2}}$. 
%\end{thm}
%The proof of this theorem directly follows from the proof of Theorem \ref{thm:sepprinc} by replacing the error system \eqref{robserragg} with \eqref{ersysos2}.

%###########################################################
\section{Example}\label{sec:ex}
\begin{figure}[t]
\begin{minipage}{\linewidth}
	\centering
	\begin{tikzpicture}[->,>=stealth',shorten >=0pt,auto,node distance=2.8cm,
	semithick,scale=1.5]
	\node (0) at(0,0.5) [circle,draw] {0};
	\node (1) at(0.5,1) [circle,draw] {1};
	\node (2) at(0.5,0) [circle,draw] {2};
	\node (3) at(1.5,1) [circle,draw] {3};
	\node (4) at(1.5,0) [circle,draw] {4};
	% Kanten
	\draw[line width=0.7pt] (0) to (1);
	\draw[line width=0.7pt] (0) to (2);
	\draw[line width=0.7pt] (1) to (3);
	\draw[line width=0.7pt] (2) to (4);
	\draw[line width=0.7pt] [bend right=10](3) to (4);
	\draw[line width=0.7pt] [bend right=10](4) to (3);
	\end{tikzpicture}
%	\hspace{1.5cm}
%	\begin{tikzpicture}[->,>=stealth',shorten >=0pt,auto,node distance=2.8cm,
%	semithick,scale=2]
%	\node (1) at(0.5,1) [circle,draw] {1};
%	\node (2) at(0.5,0) [circle,draw] {2};
%	\node (3) at(1.5,1) [circle,draw] {3};
%	\node (4) at(1.5,0) [circle,draw] {4};
%	% Kanten
%	\draw[line width=0.7pt] [bend right=10](1) to (2);
%	\draw[line width=0.7pt] [bend right=10](2) to (1);
%	\draw[line width=0.7pt] (1) to  (3);
%	\draw[line width=0.7pt] (2) to  (4);
%	\draw[line width=0.7pt] [bend right=10](3) to (4);
%	\draw[line width=0.7pt] [bend right=10](4) to (3);
%	\end{tikzpicture}
	\caption{Communication graph $\mathcal{G}$ for cooperative output regulation. The global reference model is agent $0$, which informs agent $1$ and $2$, while agent $3$ and $4$ are uninformed agents. %Right: Communication graph $\mathcal{G}_s$ for leaderless output synchronization.
	}\label{fig1}
\end{minipage}
\end{figure}
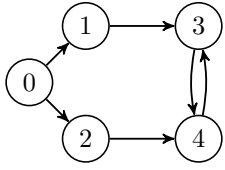

\begin{figure}[t!]
	\centering
	\input{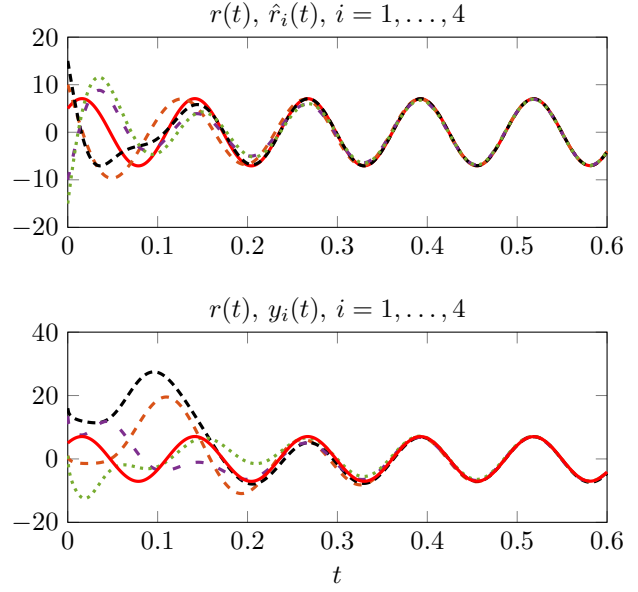}\vspace*{1ex}
	\caption{Upper plot: Estimates $\hat{r}_1$ (\ref{r.data2}), $\hat{r}_2$ (\ref{r.data3}), $\hat{r}_3$ (\ref{r.data4}) and $\hat{r}_4$(\ref{r.data5}) of the reference observers $r$ (\ref{r.data1}). Lower plot: Tracking behaviour of the agents output $y_1$ (\ref{yl.data2}), $y_2$ (\ref{yl.data3}), $y_3$ (\ref{yl.data4}) and $y_4$ (\ref{yl.data5}) for $r$ (\ref{yl.data6}) in the presence of local sinusoidal disturbances.}\label{fig:cor}
\end{figure}

The results of the paper are demonstrated for a MAS consisting of $N = 4$ parabolic agents. In what follows the \emph{normalized bump function} $\theta_\omega(t) = \e^{-((1-t)t)^{-\omega}}/\e^{-(\frac{1}{4})^{-\omega}}$ and the \emph{smooth step function} $\Theta_\omega(t) = \int_{0}^{t}\theta_\omega(\tau)\mathrm{d}\tau/\int_{0}^{1}\theta_\omega(\tau)\mathrm{d}\tau$, $t \in (0,1)$, with $\Theta_{\omega}(t) = 0$, $t = 0$, and $\Theta_{\omega}(t) = 1$, $t \geq 1$, are utilized as Gevrey functions of the order $\alpha = 1 + \tfrac{1}{\omega}$ (see \cite[App. B.1]{Me13}). The agents \textcolor{blue}{have the spatial domains $\Omega_i \in \{[0,0.9],[0,0.7],[0,0.8],[0,1]\}$, in which $\Omega_i$ is the $i$-th element of the set.} Their parameters are \textcolor{blue}{$\lambda_i(z) \in  \{0.81+0.9z,0.5+0.5z,2+z,3+\sin z\}$, $\phi_i(z) \in  \{z,0.7z+1,z,z^2\}$}, $a_i(z,t) \in  \{0.5\sin(2\pi(1.11z+4t))-28, \theta_{1.3}(t)+z-25, \sin(2\pi(z+5t))-26, \theta_{1.2}(t)+\sin z-30\}$, $q_i(t) \in \{1-\sin t,\theta_{1.2}(t)+1,3-\sin t,\theta_{1.4}(t)\}$ and \textcolor{blue}{$q^{\ell}_i(t) \in \{\sin t,\sin t + 1,\sin(t^2),1+\sin(\pi t)\}$}. The disturbance input locations of the agents are characterized by $g_{i,1} \in \{1+z+\sin3t,2\!+\!\sin4t,3\!+\!\sin2t,2\sin 2t\}$, $g_{i,2} \in \{2\!+\!2z\!+\!\sin5t,4\!+\!\sin6t,6\!+\!\sin3t,1\!+\!\sin t\}$, $g_{i,3} = \{\Theta_{1.3}\!+\!2z\!+\!2,\Theta_{1.6}\!+\!4,\Theta_{1.1}\!+\!6,\Theta_{1.3}\!+\!1\}$, $g_{i,4} \in \{\Theta_{1.3}\!+\!\sin z+2,\Theta_{1.4}\!+\!5,\Theta_{1.6}\!+\!7,\Theta_{1.5}\}$,
%\begin{equation}
% \begin{bmatrix}
%  g_{1,1}\\
%  g_{2,1}\\
%  g_{3,1}\\
%  g_{4,1} 
% \end{bmatrix}
% =
% \begin{bmatrix}
% 1\!+\!z\!+\!\sin3t\\
% 2\!+\!\sin4t\\
% 3\!+\!\sin2t\\
% 2\sin 2t
% \end{bmatrix}, \;
% \begin{bmatrix}
% g_{1,2}\\
% g_{2,2}\\
% g_{3,2}\\
% g_{4,2} 
% \end{bmatrix}
% =
% \begin{bmatrix}
% 2\!+\!2z\!+\!\sin5t\\
% 4\!+\!\sin6t\\
% 6\!+\!\sin3t\\
% 1\!+\!\sin t
% \end{bmatrix}  
% \end{equation}
%and
%\begin{equation}
%\begin{bmatrix}
%g_{1,3}\\
%g_{2,3}\\
%g_{3,3}\\
%g_{4,3} 
%\end{bmatrix}
%=
%\begin{bmatrix}
%\Theta_{1.3}\!+\!2z\!+\!2\\
%\Theta_{1.6}\!+\!4\\
%\Theta_{1.1}\!+\!6\\
%\Theta_{1.3}\!+\!1
%\end{bmatrix}, \;
%\begin{bmatrix}
%g_{1,4}\\
%g_{2,4}\\
%g_{3,4}\\
%g_{4,4} 
%\end{bmatrix}
%=
%\begin{bmatrix}
%\Theta_{1.3}\!+\!\sin z\!+\!2\\
%\Theta_{1.4}\!+\!5\\
%\Theta_{1.6}\!+\!7\\
%\Theta_{1.5}
%\end{bmatrix}, 
%\end{equation} 
in which for notational convenience the time arguments are omitted. Obviously, the Gevrey order of the system parameters ranges in the interval $[1,2)$.
%\textbf{Cooperative output regulation problem.} 
The global reference model (agent $0$) is given by
\begin{equation}
 S = \begin{bmatrix}
 0  & -50\\
 50 & 0
 \end{bmatrix}, \quad
 p\t = \begin{bmatrix}
 1 & 0
 \end{bmatrix},
\end{equation}
which generates sinusoidal reference inputs $r(t) = r_0\sin(50 t + \varphi_r)$, $r_0, \varphi_r \in \mathbb{R}$. The local sinusoidal disturbances $d_i(t) = d_{i0}\sin(10t + \varphi_i)$, $d_{i0}, \varphi_i \in \mathbb{R}$, $i = 1,\ldots,4$, are generated by \eqref{smodel}, where $S_{i} = S_d$ with
\begin{equation}
 S_d = 
 \begin{bmatrix}
  0  & -10\\
  10 &   0
 \end{bmatrix},\;
 \begin{bmatrix}
  p_1\t\\
  p_2\t\\
 \end{bmatrix} 
 =
 \begin{bmatrix}
 1   & 2\\
 2   & 1\\
 \end{bmatrix},\;
 \begin{bmatrix}
 p_3\t\\
 p_4\t
 \end{bmatrix} 
 =
 \begin{bmatrix}
 1.5 & 1\\
 2   & 2
\end{bmatrix}.
\end{equation}
The communication of the agents is described by the digraph shown in Figure \ref{fig1}. 
The Laplacian matrix associated with the digraph $\mathcal{G}$
%is
%\begin{equation}
%L_{\mathcal{G}} = 
%\begin{bmatrix}
% 0 & 0& 0& 0& 0\\ 
%-1 & 1& 0& 0& 0\\ 
%-2 & 0& 2& 0& 0\\ 
% 0 & -1& 0& 2& -1\\ 
% 0 & 0& -1& -1& 2
%\end{bmatrix},
%\end{equation}
has the spectrum $\sigma(L_{\mathcal{G}}) =\{3,1,1,2,0\}$ verifying that $\mathcal{G}$ is connected (see, e.\:g., \cite[Th. 2.1]{Lu19}). The observer gains $l_1$ and $l_2$ for the local reference observers assign the eigenvalue sets $\sigma(S_r-  l_1a_{10}p\t) = \{-10\pm\text{j}48.99\}$ and $\sigma(S_r - l_2a_{20}p\t) = \{-20\pm\text{j}60\}$. For the cooperative reference observer the common observer gain $l_w$ follows from solving the Riccati equation of Theorem \ref{theo:robsco} for $a = 500$ and $\nu = 1$. After choosing the design parameter $\mu_1 = 28$, $\mu_2 = 26$, $\mu_3 = 26$ and  $\mu_4 = 30$ the kernel equations \eqref{keq} are solved with the method of successive approximations in MATLAB (see \cite{Vaz08,Me09,Me13}). For this a relative error of $10^{-8}$ was utilized as stopping criterion. Then, the local state feedback regulators follow from the solution of the regulator equations \eqref{greq} and \eqref{lreq}, in which $15$ iteration were taken into account to solve \eqref{lreq}. The kernel equations \eqref{obsker} are solved for $\bar{\mu}_1 = 30$, $\bar{\mu}_2 = 30$, $\bar{\mu}_3 = 30$ and $\bar{\mu}_4 = 35$. Then, the solution of the IBVP \eqref{decoupleq} are calculated numerically. For this, the corresponding ICs are determined by simulating \eqref{decoupleq} for $\gamma_i(z,0) = [1 \;\; 1]\t$, $i = 1,\ldots,4$, by making use of a FEM model and utilizing the resulting $\gamma_i(z,1)$ to determine \eqref{decoupltraf}. In order to evaluate the uniform observability of the pairs $(\gamma_i\t(1,t),S_i)$, the determinants $\det Q_{o,i}(t)$ are plotted and 
%\begin{equation}
% \det Q_{o,i}(t)
% = 
% \det\begin{bmatrix}
%  \gamma_i\t(1,t)\\
%  \mathcal{M}_{S_i}\gamma_i\t(1,t)
% \end{bmatrix}, \quad i = 1,\ldots,4
%\end{equation}
%are plotted in Figure \ref{fig:obs},
it is checked that they are unequal to zero, which verifies uniform observability on the time interval $[0,1]$. The corresponding observer gains $l_{v_i}(t)$, $i = 1,\ldots,4$, are computed by time discretization and assigning the eigenvalue set $\sigma(S_i - l_{v_i}(t_k)\gamma_i\t(1,t_k)) = \{-30,-30\}$ at each time instant $t_k = kT$, $T = 1/150$, $k \in \mathbb{N}_0$. Then, the observer gains are obtained from a linear interpolation. This simple approach is justified by the fact that the assigned dynamics is much faster than the dynamics of the time-varying parameters. The stability of the resulting observer error dynamics is then validated in simulations.
%\begin{figure}[t!]
%	\centering
%	\input{plotZ/detQ}
%	\caption{Plots of $\det Q_{o,1}(t)$ (\ref{detq.data2}), $\det Q_{o,2}(t)$ (\ref{detq.data3}), $\det Q_{o,3}(t)$ (\ref{detq.data4}) and $\det Q_{o,4}(t)$ (\ref{detq.data5}) for verifying the uniform observability of the pairs $(\gamma_i\t(1,t),S_i)$, $i = 1,\ldots,4$.}\label{fig:obs}
%\end{figure}\vspace*{1ex}
The resulting networked controlled system is simulated for the agent IC $x_{i}(z,0) = 1$, the disturbance observer IC $\hat{x}_i(z,0) = 0$, $\hat{v}_i(0) = 0$, $i = 1,\ldots,4$. The  reference observers have the IC $\hat{w}_{1}(0) = [10 \;\; 10]\t$, $\hat{w}_{2}(0) = [15 \;\; 15]\t$, $\hat{w}_{3}(0) = -[10 \;\; 10]\t$, $\hat{w}_{4}(0) = -[15 \;\; 15]\t$. The reference input $r$ is specified by $w(0) = [5 \;\; -5]\t$ and the local disturbances follow from the IC $v_i(0) = [5 \;\; 5]\t$, $i = 1,\ldots,4$. For the simulation of the PDEs a FEM model with $101$ grid points is applied \textcolor{blue}{and the agents are considered in original coordinates}. The upper plot in Figure  \ref{fig:cor} shows the synchronization of the local reference models with the global reference model. With this, the resulting lower plot in the same figure verifies cooperative output regulation for the considered MAS. It can be seen from the simulation results that after the local reference models have estimated the output $r$ of the global reference model at about $t = 0.4$ also the outputs of the agents exactly synchronize with $r$ in the presence of disturbances.

\section{Concluding remarks}
There are no obstacles to extend the presented results to other types of distributed-parameter agents such as hyperbolic or mixed parabolic/hyperbolic systems. Then, only the local regulator design has to be changed. Furthermore, by omitting the leader also the \emph{leaderless output synchronization problem} can readily be solved by extending the presented approach. If the outputs to be controlled are available for measurement, then local output feedback regulators on the basis of the internal model principle can by utilized to ensure robust cooperative output regulation. The corresponding regulator design for the considered class of time-varying DPS, however, is impeded by determining a suitable internal model.

\bibliographystyle{plain}        % Include this if you use bibtex 
\bibliography{mybib}           % and a bib file to produce the 
                                 % bibliography (preferred). The
                                 % correct style is generated by
                                 % Elsevier at the time of printing.

%\begin{thebibliography}{99}     % Otherwise use the  
                                 % thebibliography environment.
                                 % Insert the full references here.
                                 % See a recent issue of Automatica 
                                 % for the style.
%  \bibitem[Heritage, 1992]{Heritage:92}
%     (1992) {\it The American Heritage. 
%     Dictionary of the American Language.}
%     Houghton Mifflin Company.
%  \bibitem[Able, 1956]{Abl:56}
%     B.~C.~Able (1956). Nucleic acid content of macroscope. 
%     {\it Nature 2}, 7--9. 
%  \bibitem[Able {\em et al.}, 1954]{AbTaRu:54}   
%     B.~C. Able, R.~A. Tagg, and M.~Rush (1954).
%     Enzyme-catalyzed cellular transanimations.
%     In A.~F.~Round, editor, 
%     {\it Advances in Enzymology Vol. 2} (125--247). 
%     New York, Academic Press.
%  \bibitem[R.~Keohane, 1958]{Keo:58}
%     R.~Keohane (1958).
%     {\it Power and Interdependence: 
%     World Politics in Transition.}
%     Boston, Little, Brown \& Co.
%  \bibitem[Powers, 1985]{Pow:85}
%     T.~Powers (1985).
%     Is there a way out?
%     {\it Harpers, June 1985}, 35--47.

%\end{thebibliography}

\appendix

\section{Proofs of lemmas and theorems}\label{app}
%\section{Proof of Lemma \ref{lem:decoupleq}}\label{sec:app}
\textbf{Proof of Lemma \ref{lem:decoupleq}.} Since the matrix $S_i \in \mathbb{R}^{n_{v_i} \times n_{v_i}}$, $i = 1,\ldots,N$, is diagonalizable by assumption, there exist $n_{v_i}$ linearly independent eigenvectors $\upsilon^i_j$ w.r.t. the eigenvalue $\nu^i_j$, $j = 1,2,\ldots,n_{v_i}$,  of $S_i$. Postmultiply \eqref{decoupleq} by the eigenvectors and define $\gamma^{i}_j = \gamma_i\t \upsilon^i_j$, $\bar{h}^{i}_{j} = -\bar{h}_{1,i}\t\upsilon^i_j$ as well as $\tilde{g}^{i}_{k,j} = \tilde{g}_{k,i}\t\upsilon_j^{i}$, $k = 2,3$.
% leads to
%\begin{subequations}\label{decoupleq2}
%	\begin{align}
%	\gamma^i_{j,t}(z,t) &= \gamma^i_{j,zz}(z,t) - (\bar{\mu}_i + \lambda^i_j)\gamma^i_j(z,t) + \bar{h}^i_j(z,t)\label{decplpdep}\\
%	\gamma^i_{j,z}(0,t) &= \tilde{g}^i_{2,j}(t), &&\hspace{-4.5cm}  t \in \mathbb{R}^+_{t_0}\label{debc1p}\\
%	\gamma^i_{j,z}(1,t) &= \tilde{g}^i_{3,j}(t), &&\hspace{-4.5cm}  t \in \mathbb{R}^+_{t_0}\label{debc2p}
%	\end{align}
%\end{subequations}
%with \eqref{decplpdep} defined on $(z,t) \in (0,1) \times \mathbb{R}^+_{t_0}$. In order to obtain homogeneous BCs,
Then, introduce the change of coordinates $\bar{\gamma}^i_j(z,t) = \gamma^i_j(z,t) - \alpha^i_j(z)\tilde{g}^i_{2,j}(t) -  \beta^i_j(z)\tilde{g}^i_{3,j}(t)$ with $\alpha^i_j(z) = z - 0.5z^2$ and $\beta^i_j(z) = 0.5z^2$. After a simple calculation shows that this leads to
\begin{subequations}\label{decoupleq2h}
	\begin{align}
	\bar{\gamma}^i_{j,t}(z,t) &= \textcolor{blue}{\lambda_i}\bar{\gamma}^i_{j,zz}(z,t) - (\bar{\mu}_i + \nu^i_j)\bar{\gamma}^i_j(z,t) + \tilde{h}^i_j(z,t)\label{decplpdeph}\\
	\bar{\gamma}^i_{j,z}(0,t) &= 0\label{debc1ph}\\
	\bar{\gamma}^i_{j,z}(1,t) &= 0\label{debc2ph}
	\end{align}
\end{subequations}
with some function $\tilde{h}^i_j \in C^2[0,1] \times G_{\alpha}(\mathbb{R}^+_{t_0})$, $1 \leq \alpha < 2$, since the elements of $\bar{h}^i_j(\cdot,t) \in C^2[0,1]$ and $\tilde{g}^i_{k,j}(t)$, $k = 2,3$, are of the same Gevrey class w.r.t. time. Therein, the fact was utilized that the time derivative of a Gevrey function does not change its Gevrey order  (see \cite[Prop. 1.4.5]{Rod93}). Define the state $\bar{\gamma}^i_{j}(t) = \{\bar{\gamma}^i_{j}(z,t), z \in [0,1]\}$, the system operator $\mathcal{A}_j^ih = \textcolor{blue}{\lambda_i}h'' - (\bar{\mu}_i + \nu^i_j)h$, $h \in D(\mathcal{A}_j^i) = \{h \in H^2(0,1)\,|\, h'(0) = h'(1) = 0\}$ and $b_j^i(t) = \{\tilde{h}^i_j(z,t), z \in [0,1]\}$. With this, \eqref{decoupleq2h} can be represented by the abstract IVP $\dot{\bar{\gamma}}^i_{j}(t) = \mathcal{A}_j^i\bar{\gamma}^i_{j}(t)  + b_j^i(t), \quad t\in \mathbb{R}_0^+, \; \bar{\gamma}^i_{j}(t_0) \in D(\mathcal{A}_j^i)$, in the state space $X_j^i = L_2(0,1)$ endowed with the usual inner product. The system operator $\mathcal{A}_j^i$ is the generator of an \emph{analytic $C_0$-semigroup} $\mathcal{T}_j^i(t)$, $t \in \overline{\mathbb{R}}^+_{t_0}$,  on $X_j^i$, because $-\mathcal{A}_j^i$ is a \emph{Sturm-Liouville operator} (see \cite{Deu15a}). Consequently, by taking  $b_j^i(t) \in C^1[0,1]$ into account the abstract IVP in question is well-posed in view of \cite[Th. 3.1.1]{Cu95} so that a unique solution exists in $L_2(0,1)$. In order to verify that $\bar{\gamma}^i_{j}(1,t)$ is of Gevrey order $\alpha$ with $1 \leq \alpha < 2$, define the operator $\mathcal{C}h = h(1)$, $D(\mathcal{C}) = D(\mathcal{A}_j^i)$. It is shown in \cite{Deu15a} that the composite operator $\mathcal{C}\mathcal{T}_j^i(t)$ is bounded on $X_j^i$. Hence, the boundary value $\bar{\gamma}^i_{j}(1,t)$ can be represented by  $\bar{\gamma}^i_{j}(1,t) = \mathcal{C}\mathcal{T}_j^i(t)\bar{\gamma}^i_{j}(t_0) + \int_{t_0}^t\mathcal{C}\mathcal{T}_j^i(t-\tau)b_j^i(\tau)\d\tau$, $t \in \overline{\mathbb{R}}_0^+$ (see \cite[Ch. 3.1]{Cu95}). Owing to the fact that $\mathcal{T}_j^i(t)$ is an analytic $C_0$-semigroup, the function $\mathcal{C}\mathcal{T}_j^i(t)h$, $h \in X_j^i$, is of Gevrey order $\alpha = 1$. Since the addition, multiplication and integration of Gevrey functions does not change the Gevrey order (see \cite[Prop. 1.4.5]{Rod93}) and $b_j^i \in G_{\alpha}(\mathbb{R}^+_{t_0})$, $1 \leq \alpha < 2$, the boundary value $\bar{\gamma}^i_{j}(1,t)$ is a Gevrey function of order $1 \leq \alpha < 2$. Then, the elements $\gamma^i_{j}(1,t)$ of $\gamma_{i}\t(1,t)$ resulting from solving the introduced change of coordinates, have the same Gevrey order, since $\tilde{g}^i_{2,j}$ and $\tilde{g}^i_{3,j}$ share this property.

%\section{Proof of Theorem \ref{th:dobs}}\label{sec:app2a}
\textbf{Proof of Theorem \ref{th:dobs}.} In order to verify exponential stability for \eqref{tfirst2}, introduce the system operator $\mathcal{A}_ih = \textcolor{blue}{\lambda_i}h'' - \bar{\mu}_ih$ with $h \in D(\mathcal{A}_i) = \{h \in H^2(0,1)\;|\; h'(0) = h'(1) = 0\}$ for \eqref{errs2t2}--\eqref{errs4t2}. Since $-\mathcal{A}_i$ is a Sturm-Liouville operator, it is the infinitesimal generator of an analytic $C_0$-semigroup $\mathcal{T}_i(t)$ on the state space $X_i =L_2(0,1)$ with the usual inner product (see \cite{Deu15a}). Hence, by defining the state $\varepsilon_{x_i}(t) = \{\varepsilon_{x_i}(z,t), z\in [0,1]\}$, the PDE subsystem \eqref{errs2t2}--\eqref{errs4t2} can be represented by the well-posed abstract IVP $\dot{\varepsilon}_{x_i}(t) = \mathcal{A}_i\varepsilon_{x_i}(t)$, $t \in \mathbb{R}^+_{t_0}$, $\varepsilon_{x_i}(t_0) \in D(\mathcal{A}_i)$, in $X_i$. Define the operator $\mathcal{C}h = h(1)$, $D(\mathcal{C}) = D(\mathcal{A}_i)$ so that $\varepsilon_{x_i}(1,t) = \mathcal{C}\mathcal{T}_i(t)\varepsilon_{x_i}(t_0)$, $t \in \overline{\mathbb{R}}^+_{t_0}$, is well-defined. In particular, $\mathcal{C}\mathcal{T}_i(t): X_i \to \mathbb{C}$ is a bounded operator, because the semigroup $\mathcal{T}_i(t)$ is analytic and $\mathcal{C}$ is relatively compact (see \cite{Deu15a}). With this, the solution of \eqref{everrt2} reads $e_{v_i}(t) \! = \! \Psi_{i}(t,t_0)e_{v_i}(t_0) \!-\! \tint_{t_0}^t\Psi_i(t,\tau)l_{v_i}(\tau)\mathcal{C}\mathcal{T}_i(\tau)\d\tau\varepsilon_{x_i}(t_0)$, in which the fundamental matrix $\Psi_i(t,\tau)$ is defined in Lemma \ref{lem:declode}.  Consider the estimates $\|\Psi_i(t,\tau)\| \leq \kappa_i\e^{-\bar{\mu}_{v_i}(t - \tau)}$ (see Lemma \ref{lem:declode}) and $|\mathcal{C}\mathcal{T}_i(\tau)| \leq d_i\e^{-\bar{\mu}_i(t-t_0)}$, $d_i > 0$, in which the latter is implied by the spectrum determined growth assumption valid for the analytic $C_0$-semigroup $\mathcal{T}_i(t)$ leading to the estimate $\|\mathcal{T}_i(t)\| \leq \bar{M}_i\e^{-\bar{\mu}_i(t-t_0)}$, $t \in \overline{\mathbb{R}}^+_{t_0}$, $\bar{M}_i \geq 1$ (see \cite{Tr75}). By taking the boundedness of $l_{v_i}(t)$ into account (see Lemma \ref{lem:declode}), one obtains the estimate $\|e_{v_i}(t)\|_{\mathbb{C}^{n_{v_i}}}  = \bar{c}_{i}\e^{-\bar{\mu}_{v_i}(t - t_0)}\|e_{v_i}(t_0)\|_{\mathbb{C}^{n_{v_i}}} + \bar{c}_{i}\tint_{t_0}^t\e^{-\bar{\mu}_{v_i}(t - \tau)}\e^{-\bar{\mu}_{i}(\tau - t_0)}\d\tau\|\varepsilon_{x_i}(t_0)\|_{L_2}$ for $e_{v_i}(t)$ and a $\bar{c}_{i} > 0$ after simple computations. By making use of the readily verifiable result $|\tint_{t_0}^t\e^{-\bar{\mu}_{v_i}(t - \tau)}\linebreak\cdot\e^{-\bar{\mu}_{i}(\tau - t_0)}\d\tau| \leq (2/|\bar{\mu}_{v_i} - \bar{\mu}_i|)\e^{-\min(\bar{\mu}_i,\bar{\mu}_{v_i})(t-t_0)}$ for $\bar{\mu}_{v_i} \neq \bar{\mu}_i$, the estimate for $\|e_{v_i}(t)\|_{\mathbb{C}^{n_{v_i}}} $ becomes 
$\|e_{v_i}(t)\|_{\mathbb{C}^{n_{v_i}}}  \leq M_i\e^{-\min(\bar{\mu}_i,\bar{\mu}_{v_i})(t-t_0)}\|\varepsilon_{o,i}(t_0)\|_{L_2}$ for $t \in \overline{\mathbb{R}}^+_{t_0}$ with $\varepsilon_{o,i} = \col{\varepsilon_{x_i},e_{v_i}}$, $\|\cdot\| = (\|\cdot\|^2_{L_2} +\linebreak \|\cdot\|^2_{\mathbb{C}^{n_{v_i}}})^{1/2}$ and an $M_i > 0$, which follows from a straightforward calculation in view of $\|e_{v_i}(t_0)\|_{\mathbb{C}^{n_{v_i}}}  \leq \|\varepsilon_{o,i}(t_0)\|$ and $\|\varepsilon_{x_i}(t_0)\|_{L_2} \leq \|\varepsilon_{o,i}(t_0)\|$. Collecting these results one directly obtains $\|\varepsilon_{o,i}(t)\| = \big(\|\varepsilon_{x_i}(t)\|_{L_2}^2 + \|e_{v_i}(t)\|_{\mathbb{C}^{n_{v_i}}}^2\big)^{\frac{1}{2}}
\leq  \|\varepsilon_{x_i}(t)\|_{L_2} + \|e_{v_i}(t)\|_{\mathbb{C}^{n_{v_i}}}
\leq \linebreak\widetilde{M}_i\e^{-\min(\bar{\mu}_i,\bar{\mu}_{v_i})(t-t_0)}\|\varepsilon_{o,i}(t_0)\|$
for $t \in \overline{\mathbb{R}}^+_{t_0}$, all $\varepsilon_{x_i}(t_0) \in D(\mathcal{A}_i)$ and $e_{v_i}(t) \in \mathbb{R}^{n_{v_i}}$. Since this result is valid for all $t_0 \geq 0$ with the same $\widetilde{M}_i$, $\bar{\mu}_i$ and $\bar{\mu}_{v_i}$, this shows uniform exponential stability. Hence, by taking the bounded invertibility of the chain of transformations into account, the same stability result holds in the original coordinates.	

%\section{Proof of Theorem \ref{thm:sepprinc}}\label{sec:app2}
\textbf{Proof of Theorem \ref{thm:sepprinc}.} Consider the state observer \eqref{dobsss}--\eqref{dobsu} for the $i$-th agent with the feedback given in the theorem and apply the transformations
\begin{subequations}\label{app:tseq}
\begin{align}
\tilde{x}_i(z,t) &= \hat{x}_i(z,t) - \tint_0^zk_i(z,\zeta,t)\hat{x}_i(\zeta,t)\d\zeta\label{app:t1}\\
\tilde{\varepsilon}_i(z,t) &= \tilde{x}_i(z,t) - \pi_i\t(z)\hat{w}_i(t) - \varphi_i\t(z,t)\hat{v}_i(t)\label{app:t2}\\
\bar{\varepsilon}_i(z,t) &= \tilde{\varepsilon}_i(z,t) - \tfrac{1}{2}z^2l_{1,i}(t)\textcolor{blue}{c_{m,i}}\tilde{e}_{x_i}(z,t).\label{app:ftraf}
\end{align}	 
\end{subequations}
Therein, the kernel $k_i(z,\zeta,t) \in \mathbb{R}$ in \eqref{app:t1} results from solving the kernel equations \eqref{keq} as well as the vectors $\pi_i(z) \in \mathbb{R}^{n_w}$  and $\varphi_i(z,t) \in \mathbb{R}^{n_{v_i}}$ in \eqref{app:t2} are the solutions of the regulator equations \eqref{greq} and \eqref{lreq}. The final transformation \eqref{app:ftraf} is utilized to homogenize the BCs for the system resulting from the transformations \eqref{app:t1} and \eqref{app:t2}. After straightforward but lengthy calculations, this leads with \eqref{decoupltraf} to the target system
\begin{subequations}\label{app:esys}
	\begin{align}
	 \dot{\bar{\varepsilon}}_i(z,t) &= \textcolor{blue}{\lambda_i}\bar{\varepsilon}''_i(z,t) - \mu_i\bar{\varepsilon}_i(z,t) -2zl_{1,i}(t)\textcolor{blue}{c_{m,i}}\varepsilon'_{x_i}(z,t)\nonumber\\
	 &\quad + \alpha_i(z,t)\varepsilon_{x_i}(z,t) + \beta_i(z,t)\varepsilon_{x_i}(1,t)\nonumber\\
	 &\quad + \delta_i\t(z,t)e_{v_i}(t) -  \pi_i\t(z)\underline{e}_{w_i}(t)\\
	 \bar{\varepsilon}'_i(0,t) &= 0\\
	 \bar{\varepsilon}'_i(1,t) &= 0
	\end{align}
\end{subequations}
for $i = 1,\ldots,N$, in which $\underline{e}_{w_i} = l_ia_{i0}p\t e_{w_i}$ for $i = 1,\ldots,n$ and $\underline{e}_{w_i} = -(e_i\t L_{32} \otimes l_wp\t )e_w^l - (e_i\t L_{33} \otimes l_wp\t)e_w^g$ for $i = n+1,\ldots,N$ (cf. \eqref{robserragg}) as well as some functions $\alpha_i, \beta_i \in [0,1] \times G_{\alpha}(\mathbb{R}^+_{t_0})$, $1 \leq \alpha < 2$,  and $\delta_i \in ([0,1] \times G_{\alpha}(\mathbb{R}^+_{t_0}))^{n_{v_i}}$.
%Define the operators $\bar{\mathcal{A}}_ih = h'' - \mu_ih$ with $h \in D(\bar{\mathcal{A}}_i) = \{h \in H^2(0,1)\;|\; h'(0) = h'(1) = 0\}$, $\mathcal{B}_{1,i}(t)h = -2(\cdot)l_{1,i}(t)h$, $\mathcal{B}_{2,i}(t)h = \alpha_i(\cdot,t)h$, $\mathcal{B}_{3,i}(t)h = \beta_i(\cdot,t)h$, 
%$h \in X_i$, $\mathcal{B}_{4,i}(t)h = \delta_i\t(\cdot,t)h$, $h \in \mathbb{C}^{n_{v_i}}$, $t \in \overline{\mathbb{R}}^+_{t_0}$, and $\mathcal{B}_{5,i}h = -\pi\t h$, $h \in \mathbb{C}^{n_w}$, in which $X_i = L_2(0,1)$ is a Hilbert space with the usual inner product. Furthermore, let $\bar{\mathcal{C}}h = h'$ and $\bar{\mathcal{C}}_1h = h(1)$, $h \in D(\bar{\mathcal{A}}_i)$. With this and introducing the state $\bar{\varepsilon}_{i}(t) = \{\bar{\varepsilon}_{i}(z,t), z \in [0,1]\}$ the error system \eqref{app:esys} can be represented by the abstract IVP
%\begin{align}\label{app:aivpesys}
% \dot{\bar{\varepsilon}}_{i}(t) &= \bar{\mathcal{A}}_i\bar{\varepsilon}_{i}(t) + \big(\mathcal{B}_{1,i}(t)\bar{\mathcal{C}} + \mathcal{B}_{2,i}(t) 
% + \mathcal{B}_{3,i}(t)\bar{\mathcal{C}}_1\big)\varepsilon_{x_i}(t)\nonumber\\
% & \quad + \mathcal{B}_{4,i}(t)e_{v_i}(t) +  \mathcal{B}_{5,i}\underline{e}_{w_i}(t)
%\end{align}
%for $t \in \mathbb{R}^+_{t_0}$, $\bar{\varepsilon}_{i}(t_0) \in  D(\bar{\mathcal{A}}_i) \subset X_i$. 
Therein, the IVP for $\varepsilon_{x_i}(t)$ can be found in proof of Theorem \ref{th:dobs} and $e_{v_i}(t)$ is the solution of \eqref{everrt2}. Finally,  $\underline{e}_{w_i}(t)$  is governed by \eqref{ewsys} for $i = 1,\ldots,n$ and by \eqref{ewgsys} for $i = n+1,\ldots,N$. %Let $\mathcal{T}_{i}(t)$ be the analytic $C_0$-semigroup introduced in Appendix \ref{sec:app2a} describing the time evolution of $\varepsilon_{x_i}(t)$. It can be verified that $\bar{\mathcal{C}}$ and $\bar{\mathcal{C}}_1$ are relatively bounded w.r.t. $\mathcal{A}_i$ (see Appendix \ref{sec:app2a} and \cite[Ch. IV/1, Sec. 2]{Kat95}). Then, with $\mathcal{T}_{i}(t)$ being an analytic $C_0$-semigroup, the composite operators $\bar{\mathcal{C}}\mathcal{T}_i(t)$ and $\bar{\mathcal{C}}_1\mathcal{T}_i(t)$ are bounded on $D(\mathcal{A}_i)$, which is shown in \cite{Deu15a}. 
By applying the reasoning of the proof for Theorem \ref{th:dobs}, uniform exponential convergence can be verified for \eqref{app:esys}. With this, the uniform exponential stability of \eqref{tfirst2}, the exponential stability of \eqref{robserragg} (see Theorem \ref{theo:robsco}) and going through the invertible chain of the boundedly invertible  transformations verifies the uniform exponential stability of the closed-loop system of Theorem \ref{thm:sepprinc} in the original coordinates  w.r.t. the norm $\|\cdot\|_{X_{\text{co}}}$.  For verifying cooperative output regulation, apply the sequence of transformations \eqref{app:tseq}  to the output tracking error $e_{y_i}(t)$ in \eqref{eyout}. This yields 
$e_{y_i}(t) = c_i\bar{\varepsilon}_i(0,t) -p\t e_{w_i}(t) + \tilde{g}\t_{4,i}(t)e_{v_i}(t)$ 
%$e_{y_i}(t) = \bar{\mathcal{C}}_0\bar{\varepsilon}_i(t) + \mathcal{H}_1e_{w_i}(t) + \mathcal{H}_{2,i}(t)e_{v_i}(t)$ with the operators  $\bar{\mathcal{C}}_0h = h(0)$, $h \in D(\bar{\mathcal{A}}_i)$,  $\mathcal{H}_1h = -p\t h$, $h \in \mathbb{R}^{n_w}$ and  $\mathcal{H}_{2,i}(t)h =  \tilde{g}\t_{4,i}(t)h$, $h \in \mathbb{R}^{n_{v_i}}$, $t \in \overline{\mathbb{R}}^+_{t_0}$, 
in the new coordinates. %Since $\bar{\mathcal{C}}_0$ is relatively bounded w.r.t. $\bar{\mathcal{A}}_i$, which is the generator of the analytic $C_0$-semigroup $\bar{\mathcal{T}}_i(t)$, the composite operator $\bar{\mathcal{C}}_0\bar{\mathcal{T}}_i(t)$ is bounded on $D(\bar{\mathcal{A}}_i)$. This follows from the same reasoning as for the related operators appearing in \eqref{app:aivpesys}.
Then, the uniform exponential stability of the closed-loop system and using the same arguments for the unbounded output operators as in \cite{Deu15a} imply \eqref{outregdualcond2} verifying cooperative output regulation.
\end{document}